\documentclass[12pt]{article}
\usepackage{amssymb}

\usepackage{graphicx}

\newtheorem{assumption}{{\bf Assumption}}[section]
\input{tcilatex}
\begin{document}

\title{Numerical Algorithms for $1$-d Backward Stochastic Differential Equations: Convergence and Simulations\footnote{This work is supported by the National Basic Research
Program of China (973 Program), No. 2007CB814902 and No.
2007CB814906.}}

\author{Shige PENG$^{a,c}$\ \ \ Mingyu XU
$^{b,c}$\thanks{Corresponding author,
Email : xumy@amss.ac.cn} \\
{\small $^a$School of Mathematics and System Science, Shandong University, }%
\\
{\small 250100, Jinan, China}\\
{\small $^b$Institute of Applied Mathematics, Academy of Mathematics
and Systems
Science,}\\
{\small Chinese Academy of Sciences, Beijing, 100080, China.}\\
{\small $^c$Department of Financial Mathematics and Control science,
School of
Mathematical Science,} \\
{\small Fudan University, Shanghai, 200433, China.}}
\date{March 12, 2008}
\maketitle

\begin{center}
\textbf{Abstract }
\end{center}

{\small In this paper we study different algorithms for backward
stochastic differential equations (BSDE in short) basing on random
walk framework for 1-dimensional Brownian motion. Implicit and
explicit schemes for both BSDE and reflected BSDE are introduced.
Then we prove the convergence of different algorithms and present
simulation results for different types of BSDEs.}

\vspace*{.8cm} \textbf{Keywords: }{\small Backward Stochastic Differential
Equations, Reflected Stochastic Differential Equations with one barrier,
Numerical algorithm, Numerical simulation}

\vspace*{.8cm} \textbf{AMS: 60H10, 34K28}

\section{Introduction}

Non-linear backward stochastic differential equations (BSDEs in
short) were firstly introduced by Pardoux and Peng (\cite{PP1990},
1990), who proved the existence and uniqueness of the adapted
solution, under smooth square integrability assumptions on the
coefficient and the terminal condition, and when the coefficient
$g(t,\omega ,y,z)$ is Lipschitz in $(y,z)$ uniformly in $(t,\omega
)$. From then on, the theory of backward stochastic differential
equations (BSDE) has been widely and rapidly developed. And many
problems in mathematical finance can be treated as BSDEs. The
natural connection between BSDE and partial differential equations
(PDE) of parabolic and elliptic types is also important
applications. It is known that only a limited number of BSDEs can be
solved explicitly. To develop numerical methods and numerical
algorithms is very helpful, both theoretically and practically.

The solution of a BSDE is a couple of progressive measurable processes $%
(Y,Z) $, which satisfies
\begin{equation}
Y_{t}=\xi +\int_{t}^{T}g(s,Y_{s},Z_{s})ds-\int_{t}^{T}Z_{s}dB_{s},
\label{(1)}
\end{equation}
where $B$ is a Brownian motion. Here $\xi $ is terminal condition
and $g$ is a generator. From \cite{PP1990}, we know that when $\xi $
is a square integrable random variable, and $g$ satisfies Lipschitz
condition and some integrability condition, BSDE (\ref{(1)}) admits
the unique solution.

The calculation and simulation of BSDEs is essentially different from those
of SDEs (see \cite{KP}). When $g$ is linear in $y$ and $z$, we may solve the
solution of BSDE by considering its dual equation, which is a forward SDE.
However for nonlinear case of $g$, we can not find the solution explicitly.
Here we describe a software package that compute our numerical solutions for
BSDEs with a convenient user-machine interface\footnote{%
The study of simulations of BSDE has been started since 1996 in
Shandong University, Mathematical Finance Laboratory directed by
PENG Shige. First simulation was done by ZHOU Haibin, then following
his works XU Mingyu worked on this software package since her master
program(from 2000). This paper is a summary of almost all algorithms
that have been used in the package. The algorithms for reflected
BSDE with two barriers will be discussed in details in another
paper.}. This package computes solutions of BSDEs, reflected BSDEs
with one or two barriers as well as BSDEs with constraints. One for
significant advantage of this package is that users have a very
convenient interface. Any users who know the ABC of BSDE can use
this package very easily. The input-output interface was also
carefully designed.

This paper is organized as follows. In Section 2, we introduce the
discretization of BSDEs, then present implicit and explicit schemes
for numerical calculation and consider their convergence. In Section
3, we show some numerical simulations. In Section 4, we consider
reflected BSDEs with one barrier which is an It\^{o} processe, by
implicit reflected scheme, explicit reflected scheme, penalized
explicit-implicit scheme and penalized explicit scheme, then we
prove the convergence of these schemes. In Section 5, we apply
penalized schemes to BSDEs with constraint on $z$ and BSDE with
solution $y$ reflecting on a function of $z$.

We should point out that there have been many recent different
algorithms for computing solutions of BSDEs and the related results
in numerical analysis, for example \cite{BP1}, \cite{BP2},
\cite{BT}, \cite{bdm}, \cite{BDM2}, \cite{C}, \cite{CMM}, \cite
{DMP}, \cite{GLX}, \cite{MPST}, \cite{MPX}, \cite{Zh1}, \cite{Zh2},
\cite{ZZ}. In contrast to these results, our method uses very simple
method.


\section{Numerical Schemes for Standard BSDEs}


Let $(\Omega ,\mathcal{F},P)$ be a complete probability space, $%
(B_{t})_{t\geq 0}$ be a $1$-dimensional Brownian motion defined on a
fixed interval $[0,T]$. We denote by $\{\mathcal{F}_{t}\}_{0\leq
t\leq T}$ the
natural filtration generated by the Brownian motion $B$, i.e., $\mathcal{F}%
_{t}=\sigma \{B_{s};0\leq s\leq t\}$ augmented with all $P$-null sets of $%
\mathcal{F}$. We consider for a fixed $n\in \mathbf{N},$
\[
B_{t}^{n}:=\sqrt{\delta }\sum_{j=1}^{[t/\delta ]}\varepsilon _{j}^{n},\;\;%
\mbox{ for all }\;0\leq t\leq T,\;\delta =\frac{T}{n},
\]
where $\{\varepsilon _{j}^{n}\}_{j=1}^{n}$ is a $\{1,-1\}$-valued i.i.d.
sequence with $P\{\varepsilon _{j}^{n}=1\}=P\{\varepsilon _{j}^{n}=-1\}=0.5$%
, i.e., a Bernoulli sequence. We set $\mathcal{G}_{j}^{n}:=\sigma
\{\varepsilon _{1}^{n},\cdots ,\varepsilon _{j}^{n}\}$ and $t_{j}=\delta j$.

Let $g:[0,T]\times \mathbb{R}\times \mathbb{R}\rightarrow \mathbb{R}$ be a
Lipschitz function in $(y,z)$ uniformly of $t$, i.e., $g$ satisfies for a
fixed $\mu >0$
\begin{eqnarray}
|g(t,y_{1},z_{1})-g(t,y_{2},z_{2})| &\leq &\mu (|y_{1}-y_{2}|+|z_{1}-z_{2}|)
\label{Lip} \\
\forall t &\in &[0,T],\forall (y_{1},z_{1}),(y_{2},z_{2})\in \mathbb{R}%
\times \mathbb{R}.  \nonumber
\end{eqnarray}
And $g(\cdot ,0,0)\ $is square integrable.

We will approximate a pair of real-valued
$(\mathcal{F}_{t})$-progressively measurable processes $(Y,Z)$
defined on $[0,T]$ such that $E[\sup_{0\leq t\leq
T}|Y_{t}|^{2}]+E[\int_{0}^{T}|Z_{t}|^{2}dt]<\infty $, which
satisfies
\begin{equation}
-dY_{t}=g(t,Y_{t},Z_{t})dt-Z_{t}dB_{t}  \label{bsde}
\end{equation}
with given terminal condition $Y_{T}=\xi
\in\mathbf{L}^2(\mathcal{F}_T)$, where $\mathbf{L}^2(\mathcal{F}_T)$
is the space of $\mathcal{F}_T$ measurable random variable
satisfying $E\left| \xi \right| ^{2}<\infty $. It is clear that $Y$
has continuous paths. An existence and
uniqueness theorem for equation (\ref{bsde}) was established in \cite{PP1990}%
, when the generator $g$ satisfies (\ref{Lip}) and $g(\cdot ,0,0)\
$is a square integrable. In many situations we are also interested
in BSDEs of the following form:
\begin{equation}
-dY_{t}=g(t,Y_{t},Z_{t})dt+dA_{t}-Z_{t}dB_{t},\ t\in [0,T],  \label{super-g}
\end{equation}
where $(A_{t})_{t\in [0,T]}$ is an $(\mathcal{F}_{t})$-predictable RCLL
process with almost surely bounded variation such that $A_{0}=0$ and $%
E[\sup_{0\leq t\leq T}|A_{t}|^{2}]<\infty $. By the standard existence and
uniqueness theorem for solutions of BSDE, for each given $A$ and $Y_{T}=\xi
\in L^{2}(\mathcal{F}_{T})$, there exists a unique pair $(Y,Z)$ for equation
(\ref{super-g}). Here $Y$ has RCLL paths. We call the triple $(Y,Z,A)$ a $g$%
-supersolution (resp. $g$-subsolution), if $A$ is an increasing process
(resp. decreasing process). It is called a $g$-solution if $A\equiv 0$. It
is easy to check that, if both $(Y,Z,A)$ and $(Y,\bar{Z},\bar{A})$ are $g$%
-supersolutions on $[0,T]$, then $(Z,A)\equiv (\bar{Z},\bar{A})$. Thus we
often call $Y$ a $g$--super(sub)solution (or $g$-solution when $A\equiv 0$)
without specifying the related $(Z,A)$.

\subsection{Implicit and Explicit Schemes for BSDEs}
We first give an assumption for discrete terminal condition $\xi^n$.
\begin{assumption}
\label{term-num}Consider $\xi $ which is $\mathcal{F}_{T}$-measurable and $%
\xi ^{n}$ which is $\mathcal{G}_{n}^{n}$-measurable, such that
\[
E[|\xi |^{2}]+\sup_{n}E[|\xi ^{n}|^{2}]<\infty
\]
and
\[
\lim_{n\rightarrow \infty }E[|\xi -\xi ^{n}|^{2}]=0.
\]
\end{assumption}
\begin{example}Set $\xi=\Phi((B_t)_{0\leq t\leq T})$, where $\Phi:\mathbf{D}_{[0,T]}:\rightarrow
R$ and satisfies Lipschitz condition. By Donsker's theorem and
Skorokhod representation theorem, there exists a probability space,
such that $\sup_{0\leq t\leq T}\left|
B_{t}^{n}-B_{t}\right| \rightarrow 0$, as $n\rightarrow \infty $, in $%
\mathbf{L}^{2}(\mathcal{F}_{T})$, since $\varepsilon _{k}$ is in $\mathbf{L}%
^{2+\delta }$. So $\xi^n:=\Phi((B^n_t)_{0\leq t\leq T})$, with
$\xi$, satisfies Assumption \ref{term-num}.
\end{example}

The numerical solution of (\ref{bsde}) is obtained by $(Y_{t}^{n},Z_{t}^{n})%
\equiv (y_{j}^{n},z_{j}^{n})$, $t\in [j\delta ,(j+1)\delta )$,
$\delta n=T$. $(y_{j}^{n},z_{j}^{n})_{0\leq j\leq n}$ is the
solution of discrete BSDE which starts from $y_{n}^{n}=\xi ^{n}$.
Our discrete BSDE on the small interval is
\begin{equation}
y_{j}^{n}=y_{j+1}^{n}+g(t_{j},y_{j}^{n},z_{j}^{n})\delta
-z_{j}^{n}\varepsilon _{j+1}^{n}\sqrt{\delta }.\;  \label{dis-bsde}
\end{equation}
Then for given $y_{j+1}^{n}$, we want to find $\mathcal{G}_{j}^{n}$%
-measurable $(y_{j}^{n},z_{j}^{n})$. The feasibility of this scheme for
small $\delta $ is due to the following easy lemma.

\begin{lemma}
\label{num-lr} Let $y_{j+1}^{n}$ be a given $\mathcal{G}_{j+1}^{n}$%
-measurable random variable. Then, when $\delta <1/\mu $, there exists a
unique $\mathcal{G}_{j}^{n}$-measurable pair $(y_{j}^{n},z_{j}^{n})$
satisfying (\ref{dis-bsde}).
\end{lemma}

\noindent \textbf{Proof}. We set $Y_{+}=y_{j+1}^{n}|_{\varepsilon
_{j+1}^{n}=1}$ and $Y_{-}=y_{j+1}^{n}|_{\varepsilon _{j+1}^{n}=-1}$. Both $%
Y_{+}$ and $Y_{-}$ are $\mathcal{G}_{j}^{n}$-measurable. Equation (\ref
{dis-bsde}) is then equivalent to the following algebraic equation:
\begin{eqnarray*}
y_{j}^{n} &=&Y_{+}+g(t_{j},y_{j}^{n},z_{j}^{n})\delta -z_{j}^{n}\sqrt{\delta }, \\
y_{j}^{n} &=&Y_{-}+g(t_{j},y_{j}^{n},z_{j}^{n})\delta
+z_{j}^{n}\sqrt{\delta }.
\end{eqnarray*}
This is equivalent to
\begin{equation}
z_{j}^{n}=\frac{1}{2\sqrt{\delta }}(Y_{+}-Y_{-})=\frac{1}{\sqrt{\delta }}%
E[y_{j+1}^{n}\varepsilon _{j+1}|\mathcal{G}_{j}^{n}].  \label{disz}
\end{equation}
and
\begin{equation}
y_{j}^{n}-g(t_{j},y_{j}^{n},z_{j}^{n})\delta =\frac{1}{2}%
(Y_{+}+Y_{-})=E[y_{j+1}^{n}|\mathcal{G}_{j}^{n}].  \label{impli}
\end{equation}
Because $g$ is assumed to be Lipschitz, the mapping $\Theta
(y)=y-g(t_{j},y,z_{j}^{n})\delta $ is strictly monotonic: when
$\delta \mu<1$,
\[
\left\langle \Theta (y)-\Theta (y^{\prime }),y-y^{\prime }\right\rangle \geq
(1-\delta \mu )\left| y-y^{\prime }\right| ^{2}>0.
\]
So there exists a unique value $y_{j}^{n}$ satisfying (\ref{impli}). $%
\square $

This lemma shows a way to solve (\ref{dis-bsde}), and we named this
algorithm as 'implicit scheme'. In many cases, $\Theta ^{-1}$ cannot be
solved explicitly. Thus we introduce the following explicit scheme by using $%
E[y_{j+1}^{n}|\mathcal{G}_{j}^{n}]$ to approximate $y_{j}^{n}$ in $g$ of (%
\ref{dis-bsde}). We set $\bar{Y}_{T}^{n}=\bar{y}_{n}^{n}=\xi ^{n}$ and,
starting from $j=n-1$, solve in following reverse order,
\begin{equation}
\bar{y}_{j}^{n}=\bar{y}_{j+1}^{n}+g(t_{j},E[\bar{y}_{j+1}^{n}|\mathcal{G}%
_{j}^{n}],\bar{z}_{j}^{n})\delta -\bar{z}_{j}^{n}\varepsilon _{j+1}^{n}\sqrt{%
\delta }.  \label{expli}
\end{equation}
Then we get,
\begin{eqnarray*}
\overline{y}_{j}^{n} &=&E[\bar{y}_{j+1}^{n}|\mathcal{G}_{j}^{n}]+g(t_{j},E[%
\bar{y}_{j+1}^{n}|\mathcal{G}_{j}^{n}],\overline{z}_{j}^{n})\delta , \\
\overline{z}_{j}^{n} &=&\frac{1}{\sqrt{\delta }}E[\bar{y}%
_{j+1}^{n}\varepsilon^n _{j+1}|\mathcal{G}_{j}^{n}]=\frac{\bar{y}%
_{j+1}^{n}|_{\varepsilon^n _{j+1}=1}-\bar{y}_{j+1}^{n}|_{\varepsilon^n _{j+1}=-1}%
}{2\sqrt{\delta }}.
\end{eqnarray*}
This explicit scheme is useful when $g$ is not linear in $y$, for
example $g(t,y,z)=sin(y)$.
\begin{example}
In pricing option, if deposit interest $r$ and loan interest $R$ are
different, we get $g(t,y,z)=ry+\sigma \theta z+(R-r)(y-z)^-$.
\end{example}

\begin{remark}
To find $g$-super(sub)solution with an increasing process $A$ as in (\ref
{super-g}), we need to consider the discretization of $A$, setting $%
A_{0}^{n}=0$, $A_{j}^{n}:=\sum_{i=0}^{j-1}E[A_{t_{i+1}}-A_{t_{i}}|\mathcal{G}%
_{i}^{n}]$. Since $A$ is an increasing process, $A_{j}^{n}$ is also
increasing. Then instead of (\ref{dis-bsde}), we get
\[
y_{j}^{n}=y_{j+1}^{n}+g(t_{j},y_{j}^{n},z_{j}^{n})\delta
+(A_{j+1}^{n}-A_{j}^{n})-z_{j}^{n}\varepsilon _{j+1}^{n}\sqrt{\delta },
\]
where $A_{j+1}^{n}-A_{j}^{n}$ is $\mathcal{G}_{j}^{n}$-measurable. Then from
implicit scheme we get
\begin{eqnarray*}
y_{j}^{n} &=&\Theta ^{-1}(E[y_{j+1}^{n}|\mathcal{G}%
_{j}^{n}]+(A_{j+1}^{n}-A_{j}^{n})), \\
z_{j}^{n} &=&\frac{1}{\sqrt{\delta }}E[y_{j+1}^{n}\varepsilon _{j+1}^{n}|%
\mathcal{G}_{j}^{n}]=\frac{y_{j+1}^{n}|_{\varepsilon
_{j+1}^{n}=1}-y_{j+1}^{n}|_{\varepsilon _{j+1}^{n}=-1}}{2\sqrt{\delta }}.
\end{eqnarray*}
And from explicit scheme, we get
\begin{eqnarray*}
\overline{y}_{j}^{n} &=&E[\bar{y}_{j+1}^{n}|\mathcal{G}_{j}^{n}]+g(t_{j},E[%
\bar{y}_{j+1}^{n}|\mathcal{G}_{j}^{n}],\overline{z}_{j}^{n})\delta
+(A_{j+1}^{n}-A_{j}^{n}), \\
\overline{z}_{j}^{n} &=&\frac{1}{\sqrt{\delta }}E[\bar{y}%
_{j+1}^{n}\varepsilon _{j+1}^{n}|\mathcal{G}_{j}^{n}]=\frac{\bar{y}%
_{j+1}^{n}|_{\varepsilon _{j+1}^{n}=1}-\bar{y}_{j+1}^{n}|_{\varepsilon
_{j+1}^{n}=-1}}{2\sqrt{\delta }}.
\end{eqnarray*}
In this paper, we will not make special efforts to study the convergence of
discrete $g$-super(sub)solution. Indeed, if we set
\[
\widetilde{y}_{j}^{n}=y_{j}^{n}+A_{j}^{n},\widetilde{z}_{j}^{n}=z_{j}^{n},\;%
\;0\leq j\leq n,
\]
then $(\widetilde{y}^{n},\widetilde{z}^{n})$ is discrete solution of
discrete BSDE with coefficient $\widetilde{g}(t,y,z)=g(t,y-A_{t},z)$. When $%
A^{n}\rightarrow A$ in certain sense, then we can get the convergence of $%
(y^{n},z^{n})$ by $(\widetilde{y}^{n},\widetilde{z}^{n})$, which is discrete
solution of a classical BSDE.

However in many cases, the increasing process $A$ is not given, it is
associated with $(Y,Z)$ in order to keep $(Y,Z)$ to satisfying certain
condition, like reflected BSDE and constraint BSDE. We will discuss them
later in this paper.
\end{remark}

\subsection{Convergence Results for Numerical Schemes for BSDEs}

We set
\[
Y_{t}^{n}=y_{[t/\delta ]}^{n},Z_{t}^{n}=z_{[t/\delta ]}^{n},\;\;\bar{Y}%
_{t}^{n}=\bar{y}_{[t/\delta ]}^{n},\;\overline{Z}_{t}^{n}=\overline{z}%
_{[t/\delta ]}^{n},\;\;\;0\leq t\leq T,
\]
where $(y_{j}^{n},z_{j}^{n})_{0\leq j\leq n}$ and $(\overline{y}_{j}^{n},%
\overline{z}_{j}^{n})_{0\leq j\leq n}$ are discrete solutions of
(\ref {dis-bsde}) by implicit and explicit schemes, respectively.

By Donsker's theorem and Skorokhod representation theorem, there
exists a probability space, such that $\sup_{0\leq t\leq T}\left|
B_{t}^{n}-B_{t}\right| \rightarrow 0$, as $n\rightarrow \infty $, in $%
\mathbf{L}^{2}(\mathcal{F}_{T})$, since $\varepsilon _{k}$ is in $\mathbf{L}%
^{2+\delta }$. Here $\mathbf{L}^{2+\delta }$ is the space of random
variable $\phi$ satisfying $E[(\phi)^{2+\delta}]<+\infty$. Then we
have

\begin{theorem}
\label{conv1}We suppose that assumptions \ref{term-num} hold and that $g$ is
Lipschitz in $y$ and $z$. Then the discrete solutions $\{(Y^{n},Z^{n})%
\}_{n=1}^{\infty }$ under the implicit scheme and $\{(\bar{Y}^{n},\bar{Z}%
^{n})\}_{n=1}^{\infty }$ under the explicit scheme converge to the solution $%
(Y,Z)$ of (\ref{bsde}) in the following senses: as $n\rightarrow \infty $,
\begin{equation}
E[\sup_{0\leq t\leq T}\left| Y_{t}^{n}-Y_{t}\right| ^{2}+\int_{0}^{T}\left|
Z_{s}^{n}-Z_{s}\right| ^{2}ds]\rightarrow 0,  \label{conv-Im}
\end{equation}
and
\begin{equation}
E[\sup_{0\leq t\leq T}\left| \bar{Y}_{t}^{n}-Y_{t}\right|
^{2}+\int_{0}^{T}\left| \bar{Z}_{s}^{n}-Z_{s}\right|
^{2}ds]\rightarrow 0. \label{conv-exp}
\end{equation}
\end{theorem}

The convergence (\ref{conv-Im}) for this implicit scheme was obtained in
2001 by a profound result of Briand, Delyon and J. M\'{e}min \cite{bdm},
which can also be found in \cite{BDM2}. From these results, the convergence (%
\ref{conv-exp}) can be derived. Before proving (\ref{conv-exp}), we first
present following lemmas.

\begin{lemma}
\label{dis-gro}Let $a$, $b$ and $\alpha $ be positive constants, $\delta b<1$
and a sequence $(v_{j})_{j=1,\ldots n}$ of positive numbers such that, for
every $j$%
\[
v_{j}+\alpha \leq a+b\delta \sum_{i=1}^{j}v_{i}.
\]
Then
\[
\sup_{j\leq n}v_{j}+\alpha \leq ae^{bT}.
\]
\end{lemma}

This is a type of Gronwall lemma for discrete cases. The proof can be found
in \cite{MPX}, so we omit it.

\begin{lemma}
\label{est-expBSDE}We assume that $\delta $ is small enough such that $%
(1+2\mu +2\mu ^{2})\delta <1$. Then
\begin{equation}
E[\sup_{j}\left| \overline{y}_{j}^{n}\right|
^{2}+\sum_{j=0}^{n-1}\left| \overline{z}_{j}^{n}\right| ^{2}\delta]
\leq C_{\xi ^{n},g}e^{(1+2\mu +2\mu ^{2})T} \label{est-exp}
\end{equation}
where $C_{\xi ^{n},g}=(1+\delta \mu )E[\left| \xi ^{n}\right|
^{2}]+\sum_{j=0}^{n-1}g^{2}(t_{j},0,0)\delta .$
\end{lemma}

\noindent \textbf{Proof. }From explicit scheme
\[
\overline{y}_{j}^{n}=\overline{y}_{j+1}^{n}+g(t_{j},E[\overline{y}_{j+1}^{n}|%
\mathcal{G}_{j}^{n}],\overline{z}_{j}^{n})\delta -\overline{z}_{j}^{n}\sqrt{%
\delta }\varepsilon _{j+1}.
\]
We have
\begin{eqnarray}
|\bar{y}_{j}^{n}|^{2}-|\overline{y}_{j+1}^{n}|^{2} &=&-|\bar{z}%
_{j}^{n}|^{2}\delta +2[\bar{y}_{j}^{n}\cdot g(t_{j},E[\overline{y}%
_{j+1}^{n}|\mathcal{G}_{j}^{n}],\overline{z}_{j}^{n})]\delta \label{dis-ito}\\
&&-|g(t_{j},E[\overline{y}_{j+1}^{n}|\mathcal{G}_{j}^{n}],\overline{z}%
_{j}^{n})|^{2}\delta ^{2} \nonumber \\
&&-2\bar{y}_{j}^{n}\bar{z}_{j}^{n}\sqrt{%
\delta }\varepsilon _{j+1}+2\bar{z}_{j}^{n}g(t_{j},E[\overline{y}_{j+1}^{n}|\mathcal{G}_{j}^{n}],\overline{z}%
_{j}^{n})\delta\sqrt{%
\delta }\varepsilon _{j+1} \nonumber
\end{eqnarray}
Taking expectation and the sum for $j=i,\cdots ,n-1$ yields
\begin{eqnarray*}
E|\bar{y}_{i}^{n}|^{2} &\leq &E|\xi ^{n}|^{2}-\sum_{j=i}^{n-1}E|\bar{z}%
_{j}^{n}|^{2}\delta \\
&&+2\delta E\sum_{j=i}^{n-1}\{|\bar{y}_{j}^{n}|\cdot (|g(t_{j},0,0)|+\mu |E[%
\overline{y}_{j+1}^{n}|\mathcal{G}_{j}^{n}]|+\mu |\bar{z}_{j}^{n}|\}.
\end{eqnarray*}
Since the last term is dominated by
\begin{eqnarray*}
\delta E\sum_{j=i}^{n-1}\{|\bar{y}_{j}^{n}|^{2}(1 &+&\mu +2\mu
^{2})+|g(t_{j},0,0)|^{2}+\mu |E[\overline{y}_{j+1}^{n}|\mathcal{G}%
_{j}^{n}]|^{2}+\frac{1}{2}|\bar{z}_{j}^{n}|^{2}\}\  \\
\ &\leq &\delta E\sum_{j=i}^{n-1}\{|\bar{y}_{j}^{n}|^{2}(1+2\mu +2\mu
^{2})+|g(t_{j},0,0)|^{2}+\frac{1}{2}|\bar{z}_{j}^{n}|^{2}\}+\mu \delta E|\xi
^{n}|^{2},\
\end{eqnarray*}
we thus have
\begin{eqnarray*}
E|\bar{y}_{i}^{n}|^{2}+\frac{1}{2}\sum_{j=i}^{n-1}E|\bar{z}%
_{j}^{n}|^{2}\delta &\leq &\sum_{j=i}^{n-1}|g(t_{j},0,0)|^{2}\delta +(1+\mu
\delta )E|\xi ^{n}|^{2} \\
&&\ +\delta (1+2\mu +2\mu ^{2})\sum_{j=i}^{n-1}E|\bar{y}_{j}^{n}|^{2}
\end{eqnarray*}
Then by Lemma \ref{dis-gro}, we obtain
\[
\sup_i E|\bar{y}_{i}^{n}|^{2}+\frac{1}{2}\sum_{j=0}^{n-1}E|\bar{z}%
_{j}^{n}|^{2}\delta \leq C_{\xi ^{n},g}e^{(1+2\mu +2\mu ^{2})T}
\]
For (\ref{est-exp}), we recall (\ref{dis-ito}), and take the sum for
$j=i,\cdots ,n-1$ and $\sup$ over $j$, then take expectation. Notice
that
$\sum_{j=0}^{i}\bar{y}_{j}^{n}\bar{z}_{j}^{n}\sqrt{%
\delta }\varepsilon _{j+1}$ and $\sum_{j=0}^{i}g(t_{j},E[\overline{y}_{j+1}^{n}|\mathcal{G}_{j}^{n}],\overline{z}%
_{j}^{n})\bar{z}_{j}^{n}\delta\sqrt{%
\delta }\varepsilon _{j+1}$ are both martingales with respect to
$\mathcal{G}^n_i$, we apply Burkholder-Davis-Gundy inequality for
them with similar techniques as before, then get
\begin{eqnarray*}
E[\sup_i|\bar{y}_{i}^{n}|^{2} ]&\leq &cC_{\xi^n, g^n} + C_\mu \delta
\sum_{j=0}^{n-1}E|\bar{y}_{j}^{n}|^{2} 
\\&\leq & cC_{\xi^n,
g^n}+C_{\mu}T \sup_j E|\bar{y}_j^n|^2
\end{eqnarray*}
With previous results, we obtain (\ref{est-exp}). $\square $ \\[0.4cm]
\textbf{Proof of Theorem \ref{conv1}. }The convergence of
$(Y^{n},Z^{n})$ to
$(Y,Z)$ is proved in \cite{bdm}. To prove (\ref{conv-exp}), the result for $(%
\overline{Y}^{n},\overline{Z}^{n})$, it suffices to prove as $n\rightarrow
\infty $,
\begin{equation}
E[\sup_{0\leq t\leq T}\left| Y_{t}^{n}-\bar{Y}_{t}^{n}\right|
^{2}+\int_{0}^{T}\left| Z_{s}^{n}-\bar{Z}_{s}^{n}\right|
^{2}ds]\rightarrow 0.  \label{diff}
\end{equation}
From (\ref{dis-bsde}) and (\ref{expli}), we have
\begin{eqnarray}
\left| y_{i}^{n}-\overline{y}_{i}^{n}\right| ^{2} &=&\left| y_{i+1}^{n}-%
\overline{y}_{i+1}^{n}\right| ^{2}-E|z_{i}^{n}-\bar{z}_{i}^{n}|^{2}\delta \label{dis-ito2}\\
&&+2[(y_{i}^{n}-\bar{y}_{i}^{n})\cdot
(g(t_{j},y_{i}^{n},z_{i}^{n})-g(t_{j},E[\overline{y}_{i+1}^{n}|\mathcal{G}%
_{i}^{n}],\overline{z}_{i}^{n}))]\delta \nonumber\\
&&-|g(t_{j},y_{i}^{n},z_{i}^{n})-g(t_{j},E[\overline{y}_{i+1}^{n}|%
\mathcal{G}_{i}^{n}],\overline{z}_{i}^{n})|^{2}\delta ^{2}-2(y_{i}^{n}-\overline{y}_{i}^{n})(z_{i}^{n}-\bar{z}_{i}^{n})\sqrt{%
\delta }\varepsilon _{j+1}\nonumber\\
&&+2(g(t_{j},y_{i}^{n},z_{i}^{n})-g(t_{j},E[\overline{y}_{i+1}^{n}|\mathcal{G}%
_{i}^{n}],\overline{z}_{i}^{n}))(z_{i}^{n}-\bar{z}_{i}^{n})\delta\sqrt{%
\delta }\varepsilon _{j+1}.\nonumber
\end{eqnarray}
Then we take expectation and the sum over $i$ from $j$ to $n-1$. With $\xi ^{n}-\overline{%
\xi }^{n}=0$, we get
\begin{eqnarray*}
E\left| y_{j}^{n}-\overline{y}_{j}^{n}\right| ^{2} &\leq &-E[\delta
\sum_{i=j}^{n-1}\left| z_{i}^{n}-\overline{z}_{i}^{n}\right| ^{2}] \\
&&+2\sum_{i=j}^{n-1}E[(y_{i}^{n}-\overline{y}%
_{i}^{n})(g(t_{j},y_{i}^{n},z_{i}^{n})-g(t_{j},E[\overline{y}_{i+1}^{n}|%
\mathcal{G}_{i}^{n}],\overline{z}_{i}^{n}))]\delta \\
&\leq &-\frac{1}{2}E[\delta \sum_{i=j}^{n-1}\left| z_{i}^{n}-\overline{z}%
_{i}^{n}\right| ^{2}]+2\mu ^{2}\delta E[\sum_{i=j}^{n-1}|y_{i}^{n}-\overline{%
y}_{i}^{n}|^{2}] \\
&&+2\mu \delta E\sum_{i=j}^{n-1}|y_{i}^{n}-\overline{y}_{i}^{n}|\cdot
|y_{i}^{n}-E[\overline{y}_{i+1}^{n}|\mathcal{G}_{i}^{n}]|.
\end{eqnarray*}
Since $\overline{y}_{i}^{n}-E[\overline{y}_{i+1}^{n}|\mathcal{G}%
_{j}^{n}]=g(t_{j},E[\overline{y}_{j+1}^{n}|\mathcal{G}_{j}^{n}],\overline{z}%
_{j}^{n})\delta $, the last term is dominated by
\[
\delta \sum_{i=j}^{n-1}(2\mu +1)E|y_{i}^{n}-\bar{y}_{i}^{n}|^{2}+%
\sum_{i=j}^{n-1}\mu ^{2}E|g(t_{j},E[\overline{y}_{j+1}^{n}|\mathcal{G}%
_{j}^{n}],\overline{z}_{j}^{n})|^{2}\delta ^{3}.
\]
But with (\ref{est-exp}), the second term is bounded by $C\delta ^{2}$. We
thus have
\[
E\left| y_{j}^{n}-\overline{y}_{j}^{n}\right| ^{2}+\frac{\delta }{2}%
E[\sum_{i=j}^{n-1}\left| z_{i}^{n}-\overline{z}_{i}^{n}\right| ^{2}]\leq
(1+2\mu +2\mu ^{2})\delta [\sum_{i=j}^{n-1}E[|y_{i}^{n}-\overline{y}%
_{i}^{n}|^{2}]+C\delta ^{2}
\]
By Lemma \ref{dis-gro}, we get
\[
\sup_{j\leq n}E\left| y_{j}^{n}-\overline{y}_{j}^{n}\right| ^{2}\leq C\delta
^{2}e^{(2\mu +2\mu ^{2}+1)T}.
\]
Then we reconsider square of the difference between the discrete
solutions of implicit scheme and explicit scheme shown in
(\ref{dis-ito2}). This time we first take the sum and $\sup_j$, then
take expectation. Using Burkholder-Davis-Gundy inequality and
similar techniques, we get
\begin{eqnarray*}
E[\sup_j |y^n_j -\bar{y}^n_j|^2]&\leq &C_{\mu} E[\delta
\sum_{i=j}^{n-1}\left| y_{i}^{n}-\overline{y}_{i}^{n}\right| ^{2}+
\delta \sum_{i=j}^{n-1}\left|
z_{i}^{n}-\overline{z}_{i}^{n}\right| ^{2}]\\
&\leq & C_{\mu}T \sup_{j\leq n}E\left|
y_{j}^{n}-\overline{y}_{j}^{n}\right| ^{2},
\end{eqnarray*}
with previous results, (\ref{conv-exp}) follows. $%
\square $\\[0.3cm]

We now prove a more general result which will be useful in proving
convergence results for schemes of reflected BSDEs. Consider the following
BSDE
\begin{eqnarray}
-dY_{t} &=&[g_{1}(t,Y_{t},Z_{t})+g_{2}(t,Y_{t},Z_{t})]dt-Z_{t}dB_{t},
\label{BSDEfg} \\
Y_{T} &=&\xi .  \nonumber
\end{eqnarray}
Here $g_{1}$ and $g_{2}$ are both Lipschitz functions. Then we have the
following implicit--explicit scheme to only replace $y_{j}^{n}$ by $%
E[y_{j+1}^{n}|\mathcal{G}_{j}^{n}]$ in $g_{2}$,
\begin{equation}
\bar{y}_{j}^{n}=\bar{y}_{j+1}^{n}+g_{1}(t_{j},\overline{y}_{j}^{n},\bar{z}%
_{j}^{n})\delta +g_{2}(t_{j},[\bar{y}_{j+1}^{n}|\mathcal{G}_{j}^{n}],\bar{z}%
_{j}^{n})\delta -\bar{z}_{j}^{n}\varepsilon _{j+1}^{n}\sqrt{\delta },
\label{imp-exp}
\end{equation}
or, equivalently,
\begin{eqnarray*}
\overline{y}_{j}^{n} &=&E[\bar{y}_{j+1}^{n}|\mathcal{G}_{j}^{n}]+g_{1}(t_{j},%
\overline{y}_{j}^{n},\bar{z}_{j}^{n})\delta +g_{2}(t_{j},E[\bar{y}_{j+1}^{n}|%
\mathcal{G}_{j}^{n}],\overline{z}_{j}^{n})\delta , \\
\overline{z}_{j}^{n} &=&\frac{1}{\sqrt{\delta }}E[\bar{y}%
_{j+1}^{n}\varepsilon _{j+1}^{n}|\mathcal{G}_{j}^{n}]=\frac{\bar{y}%
_{j+1}^{n}|_{\varepsilon _{j+1}^{n}=1}-\bar{y}_{j+1}^{n}|_{\varepsilon
_{j+1}^{n}=-1}}{2\sqrt{\delta }}.
\end{eqnarray*}
We also set $\bar{Y}_{t}^{n}=\bar{y}_{[t/\delta ]}^{n},\;\bar{Z}_{t}^{n}=%
\bar{z}_{[t/\delta ]}^{n}$,$\;0\leq t\leq T$. Meanwhile we consider the
fully implicit scheme
\[
y_{j}^{n}=y_{j+1}^{n}+g_{1}(t_{j},y_{j}^{n},z_{j}^{n})\delta
+g_{2}(t_{j},y_{j}^{n},z_{j}^{n})\delta -z_{j}^{n}\varepsilon _{j+1}^{n}%
\sqrt{\delta },
\]
and let $Y_{t}^{n}=y_{[t/\delta ]}^{n},\;Z_{t}^{n}=z_{[t/\delta ]}^{n}$,$%
\;0\leq t\leq 1$.

\begin{proposition}
\label{ThCovIE}Under same assumptions of Theorem \ref{conv1}, assume $g_{1}$
and $g_{2}$ are Lipschitz functions. Let $(Y,Z)$ be the solution of BSDE (%
\ref{BSDEfg}). Then as $n\rightarrow \infty $,
\begin{equation}
E[\sup_{0\leq t\leq T}\left| \bar{Y}_{t}^{n}-Y_{t}\right|
^{2}+\int_{0}^{T}\left| \bar{Z}_{s}^{n}-Z_{s}\right|
^{2}ds]\rightarrow 0. \label{conv-impexp}
\end{equation}
Moreover there exists a constant $C_{2}$ depending on $T$ and $\mu _{2}$
which is Lipschitz constant of $g_{2}$, such that
\[
E[\sup_{0\leq t\leq T}\left| \bar{Y}_{t}^{n}-Y_{t}^{n}\right|
^{2}+\int_{0}^{T}\left| \bar{Z}_{s}^{n}-Z_{s}^{n}\right| ^{2}ds]\leq
C_{2}\delta ^{2}\text{.}
\]
\end{proposition}

The proof is similar to that of theorem \ref{conv1} and we omit it.
\begin{remark}
This scheme is very useful. For example, we will use it for
penalization BSDE, which will be discusses in section 4.1.
\end{remark}

\section{Simulation Results for BSDEs}


We consider the terminal condition $Y_{T}=\xi $ which is a function of $%
B_{T} $: $Y_{T}=\xi =\Phi (B_{T})$. In this case we set $y_{n}^{n}=\xi
^{n}=\Phi (B_{n\delta }^{n})$. It can be checked that our explicit schemes (%
\ref{expli}) (as well as the implicit scheme) will automatically derive
\[
y_{j}^{n}:=u(j,B_{j\delta }^{n})=u(j,\sqrt{\delta
}\sum_{i=1}^{j}\varepsilon
_{i}^{n}),\;z_{j}^{n}=v(j,B_{j\delta }^{n})=v(j,\sqrt{\delta }%
\sum_{i=1}^{j}\varepsilon _{i}^{n}).
\]
Since $B_{j\delta }^{n}$ takes on $j+1$ different values, the whole solution
$\{y_{j}^{n},z_{j}^{n}\}_{0\leq j\leq n-1}$ is a $2$--vector with $\frac{%
n\times (n+1)}{2}$ values. For convenience, we set $T=1$ in our
simulation part.

Applying the above numerical schemes, we have developed a Matlab
toolbox for calculating and simulating solutions of BSDEs. This
toolbox starts with a Matlab figure window with input area for
generator $g=g(t,y,z)$ and terminal function $\xi =\Phi (x)$, where
$x$ stands for $B_{T}$. Here $g$ and $\Phi $ can be any functions
accepted by Matlab. These toolboxes can be downloaded from
http://159.226.47.50:8080/iam/xumingyu/English.jsp, by clicking
'Preprint' on the left side.

Here we consider the case: $g(t,y,z)=-5\left| y+z\right| $, $\xi =\Phi
(B_{1})=\sin (\left| B_{1}\right| )$. After inputting these parameters of a
BSDE, the numerical calculation for the BSDE are launched after clicking the
button ``calculate''. When the toolbox indicates ``the calculation is
complete'', clicking any other button in button-area will produce different
types of simulations, i.e., clicking ``progress''\ will generate a figure
displaying the dynamic evolution of backward calculation of states $%
y_{j}^{n} $ which starts from $j=n$ and ends at $j=0$.

Clicking the button ``B.M. and solution y''\ will produce the
dynamic simulation of $(t,B_{t},Y_{t})$, shown in Figure 1.
\begin{center}
\centering\includegraphics[totalheight=100mm]{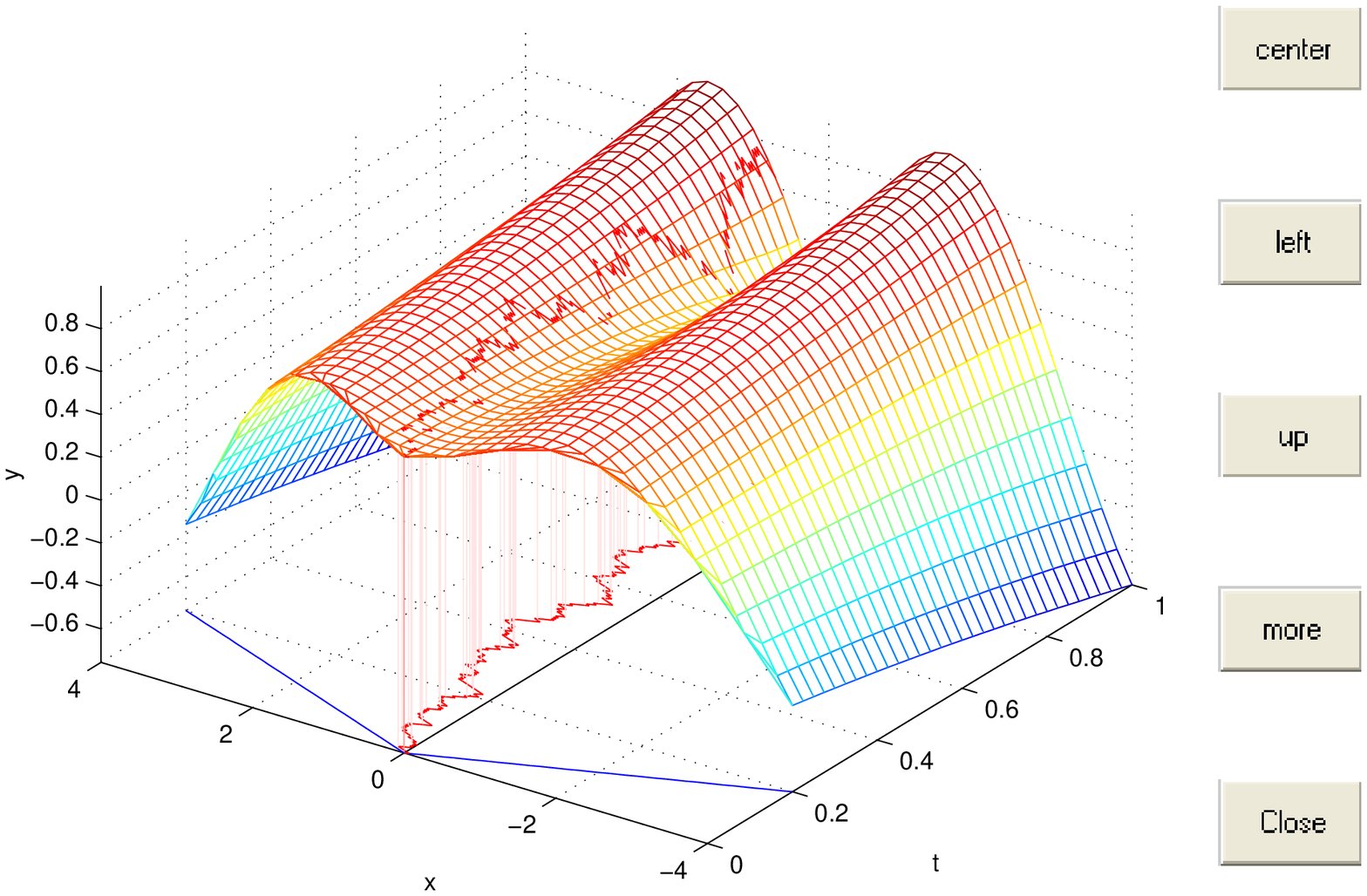} \\[0pt]
Figure 1: The solution surface with one trajectory
\end{center}
Here a trajectory of $%
Y_{t}$ runs on a colored 3-dimensional surface represented
$u=u(t,x)$, where $x$ stands for the space of Brownian motion $B$.

Clicking ``solution (y,z)''\ will generate another Matlab figure,
displayed
in Figure 2. This figure shows the 3-dimensional dynamic trajectories of $%
(t,B_{t},Y_{t})$ and $(t,B_{t},Z_{t})$ and, simultaneously, 2-dimensional
trajectories of $(t,Y_{t})$ and $(t,Z_{t})$. And there are two groups of
trajectories on the figure

\begin{center}
\centering\includegraphics[totalheight=100mm]{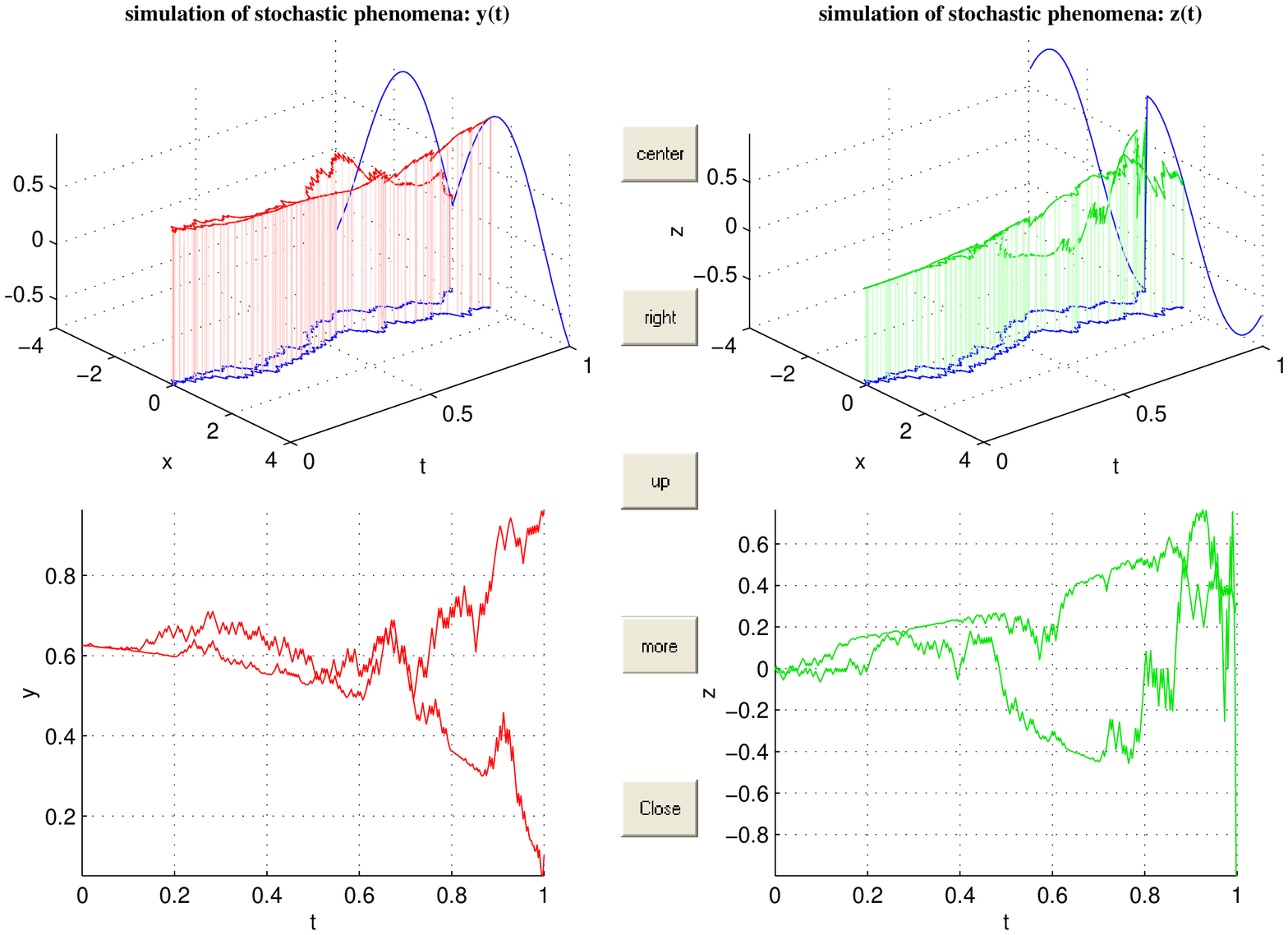} \\[0pt]
Figure 2: The trajectories of the solution
\end{center}

We now compare some numerical solutions calculated by these
algorithms: implicit scheme, explicit scheme and Monte-Carlo method
in some particular situations.

\textbf{Case I}. If $g$ is a linear function $(y,z)$:
$g(s,y,z)=by+cz+r$. The solution $Y_{0}$ of the BSDE is
\[
Y_{0}=\exp ((b-\frac{1}{2}c^{2})T)E[\xi \exp (cB_{T})]+\frac{r}{b}[\exp
(bT)-1].
\]

\begin{example} Set
$b=$ $c=$ $r=1$, $\xi =\sin (\left| B_{T}\right| )$. The numerical results
obtained with the implicit and explicit schemes are shown in the following
table:

\begin{center}
\begin{tabular}{|cccccc|}
\hline
\multicolumn{1}{|c|}{$n$} & \multicolumn{1}{c|}{100} & \multicolumn{1}{c|}{
500} & \multicolumn{1}{c|}{1000} & \multicolumn{1}{c|}{2000} & 5000 \\ \hline
\multicolumn{1}{|c|}{$Y_{0}^{n}$} & \multicolumn{1}{c|}{3.5106} &
\multicolumn{1}{c|}{3.4916} & \multicolumn{1}{c|}{3.4879} &
\multicolumn{1}{c|}{3.4866} & 3.4859 \\ \hline
\multicolumn{1}{|c|}{$\bar{Y}_{0}^{n}$} & \multicolumn{1}{c|}{3.4171} &
\multicolumn{1}{c|}{3.4716} & \multicolumn{1}{c|}{3.4785} &
\multicolumn{1}{c|}{3.4819} & 3.4840 \\ \hline
\end{tabular}
\end{center}

The exact solution is expressed by $Y_{0}=\exp (\frac{1}{2})E[\sin
(\left| B_{1}\right| )\exp (B_{1})]+\exp (1)-1$. We apply the
Monte-Carlo method, with 10,000,000 samples, to calculate $Y_{0}$.
The result is $Y_{0}=3.4850.$
\end{example}

\begin{example} Set
$b=$ $c=1$, $r=0$, $\xi =\left| B_{T}\right| $. The numerical results
obtained with the implicit and explicit schemes are:

\begin{center}
\begin{tabular}{|cccccc|}
\hline
\multicolumn{1}{|c|}{$n$} & \multicolumn{1}{c|}{100} & \multicolumn{1}{c|}{
500} & \multicolumn{1}{c|}{1000} & \multicolumn{1}{c|}{2000} & 5000 \\ \hline
\multicolumn{1}{|c|}{$Y_{0}$} & \multicolumn{1}{c|}{3.1806} &
\multicolumn{1}{c|}{3.1731} & \multicolumn{1}{c|}{3.1722} &
\multicolumn{1}{c|}{3.1719} & 3.1714 \\ \hline
\multicolumn{1}{|c|}{$\bar{Y}_{0}$} & \multicolumn{1}{c|}{3.0818} &
\multicolumn{1}{c|}{3.1531} & \multicolumn{1}{c|}{3.1621} &
\multicolumn{1}{c|}{3.1667} & 3.1694 \\ \hline
\end{tabular}
\end{center}

Applying Monte-Carlo method with 10,000,000 samples to the exact solution \\$%
Y_{0}=\exp (\frac{1}{2})E[\left| B_{1}\right| \exp (B_{1})]$, we get $%
Y_{0}=3.1710$.
\end{example}

\textbf{Case II. } If $g=\frac{1}{2}z^{2}$, then we have the exact solution $%
Y_{0}=\ln (E[\exp (\xi )])$. Since $g$ does not depend on $y$, implicit
schemes and explicit scheme give same results.

\begin{example}
For $\xi =\sin (\left| B_{1}\right| )$, applying the implicit scheme, we
obtain:

\begin{center}
\begin{tabular}{|cccccc|}
\hline
\multicolumn{1}{|c|}{$n$} & \multicolumn{1}{c|}{100} & \multicolumn{1}{c|}{
400} & \multicolumn{1}{c|}{800} & \multicolumn{1}{c|}{1000} & 2000 \\ \hline
\multicolumn{1}{|c|}{$Y_{0}^{n}$} & \multicolumn{1}{c|}{0.6249} &
\multicolumn{1}{c|}{0.6253} & \multicolumn{1}{c|}{0.6254} &
\multicolumn{1}{c|}{0.6254} & 0.6255 \\ \hline
\end{tabular}
\end{center}

By Monte-Carlo method with 10,000,000 samples to the exact expression $%
Y_{0}=\ln (E[\exp (\sin (\left| B_{1}\right| )])$, we get $Y_{0}=0.6255$.
\end{example}

\section{Reflected BSDEs}

\subsection{Algorithms for reflected BSDEs with one barrier}

In this section, we discuss the algorithms for reflected BSDEs with one
continuous lower barrier $L$. A solution of such equation is a triple $%
(Y,Z,K)$ on $[0,T]$ satisfying $E[\sup_{0\leq t\leq T}\left| Y_{t}\right|
^{2}+\int_{0}^{T}\left| Z_{s}\right| ^{2}ds+\left| K_{T}\right| ^{2}]<\infty
$ and
\begin{eqnarray}
Y_{t} &=&\xi
+\int_{t}^{T}g(s,Y_{s},Z_{s})ds+K_{T}-K_{t}-\int_{t}^{T}Z_{s}dB_{s},
\label{RBSDE1b1} \\
Y_{t} &\geq &L_{t},\;0\leq t\leq T,\text{ \ with }%
\int_{0}^{T}(Y_{t}-L_{t})dK_{t}=0.  \nonumber
\end{eqnarray}
In \cite{EKPPQ}, existence and uniqueness of the solution of this equation
is proved when $g$ satisfies Lipschitz condition  (\ref{Lip}) and $E[\left|
\xi \right| ^{2}+\int_{0}^{T}g^{2}(t,0,0)dt+\sup_{0\leq t\leq
T}(L_{t}^{+})^{2}]<\infty $. Here we consider the case when $L_{t}$ is an It%
\^{o} process, i.e. $L_{t}=L_{0}+\int_{0}^{t}l_{s}ds+\int_{0}^{t}\sigma
_{s}dB_{s}$, $0\leq t\leq T$ and $\xi =\Phi ((B_{s})_{0\leq s\leq T})$
satisfying requires of integrability, for convenience of discretization of
processes.

\begin{remark}
We call a progressively measurable process $\phi _{t}$ is in space $\mathbf{S%
}^{2}(0,T)$, if it satisfies $E[\sup_{0\leq t\leq T}\left| \phi _{t}\right|
^{2}]<\infty $. If a predictable process $\phi _{t}$ is in space $\mathbf{L}%
_{\mathcal{F}}^{2}(0,T)$, then it satisfies $E[\int_{0}^{T}\left| \phi
_{s}\right| ^{2}ds]<\infty $. And we define a space of $\mathcal{F}_{t}$%
-measurable random variables $\xi $, which satisfies $E[\left| \xi \right|
^{\beta }]<\infty $, as $\mathbf{L}^{\beta }(\mathcal{F}_{t})$, for $\beta
\in \mathbf{R}^{+}$.
\end{remark}

Following the same discretization introduced in section 2, we will
approximate the solution of reflected BSDE. On the small interval $[j\delta
,(j+1)\delta ]$, the equation (\ref{RBSDE1b1}) can be approximated by the
discrete equation
\begin{eqnarray}
y_{j}^{n} &=&y_{j+1}^{n}+g(t_{j},y_{j}^{n},z_{j}^{n})\delta
+d_{j}^{n}-z_{j}^{n}\varepsilon _{j+1}^{n}\sqrt{\delta },  \label{dis-rbsde}
\\
y_{j}^{n} &\geq
&L_{j}^{n},(y_{j}^{n}-L_{j}^{n})d_{j}^{n}=0,\nonumber
\end{eqnarray}
where $d_{j}^{n}=K_{t_{j+1}}-K_{t_{j}}$, and $L_{j}^{n}=L_{0}+\delta
\sum_{i=0}^{j-1}l_{t_{i}}+\sum_{i=0}^{j-1}\sigma _{t_{i}}\varepsilon _{i+1}^{n}%
\sqrt{\delta }$. Here (\ref{dis-rbsde}) is called discrete reflected BSDE in
\cite{MPX}, with terminal value $\xi ^{n}=\Phi ((\sum_{i=0}^{j}\varepsilon
_{i+1}^{n}\sqrt{\delta })_{0\leq j\leq n})$.

\begin{remark}
When $L_{t}=\psi (t,B_{t})$ with $\psi \in C^{1,2}([0,T]\times \mathbf{R})$,
by It\^{o} formula, we know that $L_{t}=L_{0}+\int_{0}^{t}(\frac{\partial }{%
\partial s}+\frac{1}{2}\frac{\partial ^{2}}{\partial x^{2}})\psi
(s,B_{s})ds+\int_{0}^{t}\frac{\partial }{\partial x}\psi (s,B_{s})dB_{s}$.
In fact, our algorithms are available for the case when the barrier $L$ is a
functional of Brownian motion, i.e. $L_{t}=\Psi (t,(B_{s})_{0\leq s\leq t})$%
, with its discrete version $L_{[t/\delta ]}^{n}=\Psi (t_{[t/\delta
]},(\sum_{k=0}^{i}\varepsilon _{k+1}^{n}\sqrt{\delta })_{0\leq i\leq
[t/\delta ]})$. In this section, we focus on It\^{o} process in order to
discuss the convergence of discrete solution.
\end{remark}

Suppose $y_{j+1}^{n}$ is known, we try to find $\mathcal{G}_{j}^{n}$%
-measurable $(y_{j}^{n},z_{j}^{n},d_{j}^{n})$ to satisfy
(\ref{dis-rbsde}).
Set $Y_{+}=y_{j+1}^{n}|_{\varepsilon _{j+1}^{n}=1}$ and $%
Y_{-}=y_{j+1}^{n}|_{\varepsilon _{j+1}^{n}=-1}$. From (\ref{dis-rbsde}), we
get immediately $z_{j}^{n}=\frac{1}{\sqrt{\delta}}E[y_{j+1}^{n}\varepsilon _{j+1}|\mathcal{G}%
_{j}^{n}]=\frac{1}{2\sqrt{\delta }}(Y_{+}-Y_{-})$. Substitute it
into the equation, our problem is changed to find
$(y_{j}^{n},d_{j}^{n})$ satisfying
\begin{eqnarray}
y_{j}^{n} &=&E[y_{j+1}^{n}|\mathcal{G}_{j}^{n}]+g(t_{j},y_{j}^{n},z_{j}^{n})%
\delta +d_{j}^{n},  \label{dis-rbsde1by} \\
y_{j}^{n} &\geq &L_{j}^{n},(y_{j}^{n}-L_{j}^{n})d_{j}^{n}=0.  \nonumber
\end{eqnarray}
Then we introduce two different schemes for this equation.

\paragraph{Implicit reflected scheme.}

First, we present the implicit reflected scheme which is introduces by
M\'{e}min, Peng and Xu in \cite{MPX}. If we consider the mapping $\Theta
(y):=y-(g(t_{j},y,z_{j}^{n})-g(t_{j},L_{j}^{n},z_{j}^{n}))\delta $, then for
$\delta $ small enough, we have
\[
\left\langle \Theta (y)-\Theta (y^{\prime }),y-y^{\prime }\right\rangle \geq
(1-\delta \mu )\left| y-y^{\prime }\right| ^{2}>0,
\]
i.e. $\Theta (y)$ is strictly increasing with $\Theta (L_{j}^{n})=L_{j}^{n}$%
, so
\[
\Theta ^{-1}(y)\geq L_{j}^{n}\Longleftrightarrow y\geq L_{j}^{n}.
\]
It follows
\begin{eqnarray*}
y_{j}^{n} &=&\Theta ^{-1}(E[y_{j+1}^{n}|\mathcal{G}_{j}^{n}]-g(t_{j},L_{j}^{n},z_{j}^{n})\delta +d_{j}^{n}), \\
d_{j}^{n} &=&\left( E[y_{j+1}^{n}|\mathcal{G}_{j}^{n}]+g(t_{j},L_{j}^{n},z_{j}^{n})%
\delta -L_{j}^{n}\right) ^{-}.
\end{eqnarray*}
Notice that
$E[y_{j+1}^{n}|\mathcal{G}_{j}^{n}]=\frac{1}{2}(Y_{+}+Y_{-})$, we
get the results.

\paragraph{Explicit reflected scheme}

Instead of solving the inverse of the mapping $\Theta $, we replace $%
y_{j}^{n}$ by $E[y_{j+1}^{n}|\mathcal{F}_{j}^{n}]$ on the right side of (\ref
{dis-rbsde1by}) to get an approximal solution. Then it follows
\begin{eqnarray}
\overline{y}_{j}^{n} &=&E[\overline{y}_{j+1}^{n}|\mathcal{G}%
_{j}^{n}]+g(t_{j},E[\overline{y}_{j+1}^{n}|\mathcal{G}_{j}^{n}]),\overline{z}_{j}^{n})\delta +%
\overline{d}_{j}^{n},  \label{exp-RBSDE} \\
\overline{d}_{j}^{n} &=&\left( E[\overline{y}_{j+1}^{n}|\mathcal{G}%
_{j}^{n}]+g(t_{j},E[\overline{y}_{j+1}^{n}|\mathcal{G}_{j}^{n}],\overline{z}_{j}^{n})\delta
-L_{j}^{n}\right) ^{-}.  \nonumber
\end{eqnarray}
Substitute
$E[y_{j+1}^{n}|\mathcal{F}_{j}^{n}]=\frac{1}{2}(Y_{+}+Y_{-})$ into
it, we get the results.
\begin{remark}
Compare with the implicit reflected scheme, the explicit reflected
scheme is much easier to compile programs for simulation. For
example $g(t,y,z)=sin(y)$.
\end{remark}

Another important numerical method is via the penalization equations
of reflected BSDE. In \cite{EKPPQ}, the authors introduced the
penalization method to prove the existence of the solution. For
$p\in \mathbf{N}$, the penalization equation with respect to the
lower barrier $L$ is
\begin{equation}
Y_{t}^{p}=\xi
+\int_{t}^{T}g(s,Y_{s}^{p},Z_{s}^{p})ds+p%
\int_{t}^{T}(Y_{s}^{p}-L_{s})^{-}ds-\int_{t}^{T}Z_{s}^{p}dB_{s},
\label{PBSDE1b}
\end{equation}
thanks to the comparison theorem for BSDE, we have $Y_{t}^{p}\leq
Y_{t}^{p+1} $, for $p\in \mathbf{N}$. Denote $K_{t}^{p}=p%
\int_{0}^{t}(Y_{s}^{p}-L_{s})^{-}ds$. Then we know following results from
\cite{EKPPQ}.

\begin{theorem}
\label{p1b}There exists a positive constant $c$ independent on $p$, such
that
\[
E[\sup_{0\leq t\leq
T}|Y_{t}^{p}-Y_{t}|^{2}+\int_{0}^{T}|Z_{t}^{p}-Z_{t}|^{2}dt+\sup_{0\leq
t\leq T}|K_{t}^{p}-K_{t}|^{2}]\leq \frac{c}{\sqrt{p}}.
\]
When $p\rightarrow \infty $, we know $Y^{p}\rightarrow Y$ in $\mathbf{S}%
^{2}(0,T)$, $Z^{p}\rightarrow Z$ in $\mathbf{L}_{\mathcal{F}}^{2}(0,T)$, $%
K^{p}\rightarrow K$ in $\mathbf{S}^{2}(0,T)$.
\end{theorem}

\paragraph{Numerical Penalization scheme}

By theorem \ref{p1b}, we know that the solution of reflected BSDE can be
approximated by the solution of penalization equations (\ref{PBSDE1b}), for
some large $p$. Then on the small time $[j\delta ,(j+1)\delta ]$, we
consider the following discrete penalized BSDE
\[
y_{j}^{p,n}=y_{j+1}^{p,n}+g(t_{j},y_{j}^{p,n},z_{j}^{p,n})\delta
+p(y_{j}^{p,n}-L_{j}^{n})^{-}\delta -z_{j}^{p,n}\sqrt{\delta }\varepsilon
_{j+1}.
\]
If we have already known $(y_{j+1}^{p,n},z_{j+1}^{p,n})$, then to solve $%
(y_{j}^{p,n},z_{j}^{n,p})$ from above equation, we first get $%
z_{j}^{p,n}=\frac{%
1}{\sqrt{\delta }}E[y_{j+1}^{p,n}\varepsilon _{j+1}^{n}|\mathcal{G}_{j}^{n}]=\frac{%
1}{2\sqrt{\delta }}(Y_{+}^{p}-Y_{-}^{p})$, where $Y_{+}^{p}=y_{j+1}^{p,n}|_{%
\varepsilon _{j+1}^{n}=1},Y_{-}^{p}=y_{j+1}^{p,n}|_{\varepsilon
_{j+1}^{n}=-1}.$

Then $y_{j}^{p,n}$ satisfies following equation
\begin{equation}
y_{j}^{p,n}=E[y_{j+1}^{p,n}|\mathcal{G}%
_{j}^{n}]+g(t_{j},y_{j}^{p,n},z_{j}^{p,n})\delta
+p(y_{j}^{p,n}-L_{j}^{n})^{-}\delta .  \label{dis-pbsde1}
\end{equation}
There are two ways to find suitable $y_{j}^{p,n}$. One is penalization
implicit scheme, i.e. to solve the equation:
\[
y_{j}^{p,n}=(\Theta ^{p})^{-1}(E[y_{j+1}^{p,n}|\mathcal{G}_{j}^{n}])=(\Theta
^{p})^{-1}(\frac{1}{2}(Y_{+}^{p}+Y_{-}^{p})).
\]
Here $\Theta ^{p}$ is a mapping, $\Theta
^{p}(y)=y-(g(t_{j},y,z_{j}^{p,n})+p(y-L_{j}^{n})^{-})\delta $. Let $%
d_{j}^{p,n}=p(y_{j}^{p,n}-L_{j}^{n})^{-}\delta $.

The other is implicit-explicit scheme, we only replace $y_{j}^{p,n}$
of $g$ in (\ref {dis-pbsde1}) by
$E[y_{j+1}^{p,n}|\mathcal{F}_{j}^{n}]$. Then we get, penalization
explicit-implicit scheme, i.e.
\begin{eqnarray*}
\overline{y}_{j}^{p,n} &=&E[\overline{y}_{j+1}^{p,n}|\mathcal{G}%
_{j}^{n}]+g(t_{j},E[\overline{y}_{j+1}^{p,n}|\mathcal{G}%
_{j}^{n}],z_{j}^{p,n})\delta  \\
&&+\frac{p\delta }{1+p\delta }(E[\overline{y}_{j+1}^{p,n}|\mathcal{G}%
_{j}^{n}]+g(t_{j},E[\overline{y}_{j+1}^{p,n}|\mathcal{G}_{j}^{n}],\overline{z%
}_{j}^{p,n})\delta -L_{j}^{n})^{-}.
\end{eqnarray*}
With $E[\overline{y}_{j+1}^{p,n}|\mathcal{F}_{j}^{n}]=\frac{1}{2}(\overline{y%
}_{j+1}^{p,n}|_{\varepsilon _{j+1}^{n}=1}+\overline{y}_{j+1}^{p,n}|_{%
\varepsilon _{j+1}^{n}=-1})$, results follow easily. And we set $%
\overline{d}_{j}^{p,n}=p(\overline{y}_{j}^{p,n}-L_{j}^{n})^{-}\delta $.

\subsection{Convergence results of different schemes for Reflected BSDE with
one barrier}

We first study the penalization scheme of reflected BSDE with one lower
barrier. For penalization implicit scheme, define $Y_{t}^{p,n}=y_{[t/\delta
]}^{p,n}$, $Z_{t}^{p,n}=z_{[t/\delta ]}^{p,n}$ and $K_{t}^{p,n}=%
\sum_{m=0}^{[t/\delta ]}d_{m}^{p,n}$. By Donsker's theorem and
Skorokhod representation theorem, there exists a probability space,
such that $\sup_{0\leq t\leq T}\left|
B_{t}^{n}-B_{t}\right| \rightarrow 0$, as $n\rightarrow \infty $, in $%
\mathbf{L}^{2}(\mathcal{F}_{T})$, since $\varepsilon _{k}$ is in $\mathbf{L}%
^{2+\delta }$. For convergence of scheme, we have

\begin{proposition}
\label{conl1}Under assumptions \ref{term-num} and $g$ satisfying Lipschitz
condition. The sequence $(Y_{t}^{p,n},Z_{t}^{p,n})$ converges to $%
(Y_{t},Z_{t})$ in following sense
\begin{equation}
\lim_{p\rightarrow \infty}\lim_{n\rightarrow \infty}E[\sup_{0\leq
t\leq T}\left| Y_{t}^{p,n}-Y_{t}\right| ^{2}+\int_{0}^{T}\left|
Z_{s}^{p,n}-Z_{s}\right| ^{2}ds]= 0,\ \label{convl1-np}
\end{equation}
 and for $0\leq
t\leq T$, $K_{t}^{p,n}%
\rightarrow K_{t}$ in $\mathbf{L}^{2}(\mathcal{F}_{t})$, as$\
n\rightarrow \infty $, $p\rightarrow \infty $.
\end{proposition}

\proof%
Since
\begin{eqnarray*}
E[\sup_{0\leq t\leq T}\left| Y_{t}^{p,n}-Y_{t}\right|
^{2}+\int_{0}^{T}\left| Z_{s}^{p,n}-Z_{s}\right| ^{2}ds] &\leq
&2E[\sup_{0\leq t\leq T}\left| Y_{t}^{p,n}-Y_{t}^{p}\right|
^{2}+\int_{0}^{T}\left| Z_{s}^{p,n}-Z_{s}^{p}\right| ^{2}ds] \\
&&+2E[\sup_{0\leq t\leq T}\left| Y_{t}^{p}-Y_{t}\right|
^{2}+\int_{0}^{T}\left| Z_{s}^{p}-Z_{s}\right| ^{2}ds],
\end{eqnarray*}
by the convergence results of numerical solutions for BSDE and penalization
method for reflected BSDE, Theorem \ref{p1b}, we know (\ref{convl1-np})
hold. For the increasing processes, we have
\[
E[(K_{t}^{p,n}-K_{t})^{2}]\leq
2E[(K_{t}^{p,n}-K_{t}^{p})^{2}]+2E[(K_{t}^{p}-K_{t})^{2}].
\]
While for fixed $p$,
\begin{eqnarray*}
K_{t}^{p,n}
&=&Y_{0}^{p,n}-Y_{t}^{p,n}-\int_{0}^{t}g(s,Y_{s}^{p,n},Z_{s}^{p,n})ds+%
\int_{0}^{t}Z_{s}^{p,n}dB_{s}^{n}, \\
K_{t}^{p}
&=&Y_{0}^{p}-Y_{t}^{p}-\int_{0}^{t}g(s,Y_{s}^{p},Z_{s}^{p})ds+%
\int_{0}^{t}Z_{s}^{p}dB_{s},
\end{eqnarray*}
from \cite{BDM2} Corollary 14, we know that $\int_{0}^{\cdot
}Z_{s}^{p,n}dB_{s}^{n}$ converges to $\int_{0}^{\cdot }Z_{s}^{p}dB_{s}$ in $%
\mathbf{S}^{2}(0,T)$, as $n\rightarrow \infty $, then with Lipschitz
condition of $g$ and (\ref{convl1-np}), we get $%
E[(K_{t}^{p,n}-K_{t}^{p})^{2}]\rightarrow 0$, as $n\rightarrow \infty $.
With convergence result of penalization methos, the result follows. $\square
$

Then we consider the penalization explicit-implicit scheme, note $\overline{Y%
}_{t}^{p,n}=\overline{y}_{[t/\delta ]}^{p,n}$, $\overline{Z}_{t}^{p,n}=%
\overline{z}_{[t/\delta ]}^{p,n}$ and $\overline{K}_{t}^{p,n}=%
\sum_{m=0}^{[t/\delta ]}\overline{d}_{m}^{p,n}$, it follows that

\begin{proposition}
Under same assumptions of Proposition \ref{conl1},
$(\overline{Y}_{t}^{p,n},\overline{Z}_{t}^{p,n})$ converges to
$(Y_{t},Z_{t})$ in following sense
\[
\lim_{p\rightarrow \infty}\lim_{n\rightarrow \infty}E[\sup_{0\leq
t\leq T}\left| \overline{Y}_{t}^{p,n}-Y_{t}\right|
^{2}+\int_{0}^{T}\left| \overline{Z}_{s}^{p,n}-Z_{s}\right| ^{2}ds]=
0,
\]
with $\overline{K}_{t}^{p,n}\rightarrow K_{t}$ in $\mathbf{L}^{2}(\mathcal{F}%
_{t})$, for $0\leq t\leq T$, as $n\rightarrow \infty $%
, $p\rightarrow \infty $.
\end{proposition}

\proof%
The convergence of $(\overline{Y}_{t}^{p,n},\overline{Z}_{t}^{p,n})$ is a
direct result of Proposition \ref{ThCovIE} and (\ref{convl1-np}). We
consider the increasing process, notice that for $0\leq t\leq T$,
\[
\overline{K}_{t}^{p,n}=\overline{Y}_{0}^{p,n}-\overline{Y}%
_{t}^{p,n}-\int_{0}^{t}g(s,\overline{Y}_{s}^{p,n},\overline{Z}%
_{s}^{p,n})ds+\int_{0}^{t}\overline{Z}_{s}^{p,n}dB_{s}^{n},
\]
with $K_{t}^{p,n}=Y_{0}^{p,n}-Y_{t}^{p,n}-%
\int_{0}^{t}g(s,Y_{s}^{p,n},Z_{s}^{p,n})ds+\int_{0}^{t}Z_{s}^{p,n}dB_{s}^{n}$%
, thanks to Lipschitz condition of $g$ and the convergence of $(\overline{Y}%
^{p,n},\overline{Z}^{p,n})$, we get $E[(K_{t}^{p,n}-\overline{K}%
_{t}^{p,n})^{2}]\rightarrow 0$, as $n\rightarrow \infty $, $p\rightarrow
\infty $. With convergence results of penalization method, results follow. $%
\square $

Now we study the convergence of reflected schemes. First for the implicit
reflected scheme, denote $Y_{t}^{n}=y_{[t/\delta ]}^{n}$, $%
Z_{t}^{n}=z_{[t/\delta ]}^{n}$, $K_{t}^{n}=\sum_{j=0}^{[t/\delta ]}d_{j}^{n}$%
, for $0\leq t\leq T$, from the results in \cite{MPX}, we know

\begin{theorem}[Theorem 3.2 in \cite{MPX}]
\label{conv-imp}Under assumption \ref{term-num} and (\ref{Lip}) for $g$, as $%
n\rightarrow +\infty ,$%
\[
E[\sup_{0\leq t\leq
T}|Y_{t}^{n}-Y_{t}|^{2}]+E\int_{0}^{T}|Z_{t}^{n}-Z_{t}|^{2}dt\rightarrow 0.
\]
\end{theorem}

For the increasing process, we have

\begin{proposition}
For $t\in [0,T]$, $E[(K_{t}-K_{t}^{n})^{2}]\rightarrow 0$, as $n\rightarrow
\infty $.
\end{proposition}

\proof%
For $t\in [0,T]$, we have
\[
E[(K_{t}-K_{t}^{n})^{2}]\leq
3E[(K_{t}-K_{t}^{p})^{2}]+3E[(K_{t}^{p}-K_{t}^{p,n})^{2}]+3E[(K_{t}^{p,n}-K_{t}^{n})^{2}]
\]
where $K^{p}$ is from penalization equation (\ref{p1b}), and $K_{t}^{p,n}$
is discrete solution of (\ref{dis-pbsde1}), with $K_{t}^{p,n}=%
\sum_{j=0}^{[t/\delta ]}d_{j}^{p,n}$. Similar as Lemma 2.5 in \cite{MPX}, we
have $E[\sup_{t}|Y_{t}^{p,n}-Y_{t}^{n}|^{2}]+E%
\int_{0}^{T}|Z_{t}^{p,n}-Z_{t}^{n}|^{2}dt\leq \frac{C_{\xi ^{n},g,L}^{R}}{%
\sqrt{p}}$, where $C_{\xi ^{n},g,L}^{R}$ only depends on $\xi ^{n}$, $g$, $L$
and $\mu $. Since
\begin{eqnarray*}
K_{t}^{p,n}
&=&Y_{0}^{p,n}-Y_{t}^{p,n}-\int_{0}^{t}g(s,Y_{s}^{p,n},Z_{s}^{p,n})ds+%
\int_{0}^{t}Z_{s}^{p,n}dB_{s}^{n}, \\
K_{t}^{n}
&=&Y_{0}^{n}-Y_{t}^{n}-\int_{0}^{t}g(s,Y_{s}^{n},Z_{s}^{n})ds+%
\int_{0}^{t}Z_{s}^{n}dB_{s}^{n},
\end{eqnarray*}
with Lipschitz condition of $g$, we deduce that $%
E[(K_{t}^{p,n}-K_{t}^{n})^{2}]\leq \frac{C_{\xi ^{n},g,L}^{R}}{\sqrt{p}}$.
It follows
\[
E[(K_{t}-K_{t}^{n})^{2}]\leq (C_{\xi ^{n},g,L}^{R}+C_{\xi ,g,L}^{R})\frac{1}{%
\sqrt{p}}+3E[(K_{t}^{p}-K_{t}^{p,n})^{2}].
\]
Since $K_{t}^{p,n}\rightarrow K_{t}^{p}$ in $\mathbf{L}^{2}(\mathcal{F}_{t})$
as $n\rightarrow \infty $, for fixed $p$, we can choose $n$ large enough to
get right side very small. Then result of $K^{n}$ follows. $\square $

Then we consider the convergence of the reflected explicit scheme. We set
\[
\;\bar{Y}_{t}^{n}=\bar{y}_{[t/\delta ]}^{n},\;\;\overline{Z}_{t}^{n}=%
\overline{z}_{[t/\delta ]}^{n},\;\;\overline{K}_{t}^{n}=\sum_{j=0}^{[t/%
\delta ]}\overline{d}_{j}^{n}\;\;\;0\leq t\leq T.
\]
First as Lemma \ref{est-expBSDE}, we have similar estimation of $\overline{y}%
_{j}$ of reflected BSDE, given by (\ref{exp-RBSDE}).

\begin{lemma}
\label{est-expRBSDE}We assume that $\delta $ is small enough such that $(%
2+2\mu +6\mu ^{2})\delta <1$. Then
\[
E[\sup_{j}\left| \overline{y}_{j}^{n}\right|
^{2}+\sum_{j=0}^{n-1}\left|
\overline{z}_{j}^{n}\right| ^{2}\delta] +E[(\sum_{j=0}^{n-1}\overline{d}%
_{j}^{n})^{2}]\leq C_{\xi ^{n},g,L}^{R}
\]
where $C_{\xi ^{n},g,L}^{R}$ only depends on $\mu $, $\ E[\left| \xi
^{n}\right| ^{2}]$, $\sum_{j=0}^{n-1}g^{2}(t_{j},0,0)\delta $ and $%
E[\sup_{j}((L_{j}^{n})^{+})^{2}]$.
\end{lemma}

\noindent\textbf{Proof.} Recall that for $j=0,1,...n-1$, $(\overline{y}_{j}^{n},\bar{z%
}_{j}^{n})$ satisfies
\begin{eqnarray}
\overline{y}_{j}^{n} &=&\overline{y}_{j+1}^{n}+g(t_{j},(E[\overline{y}%
_{j+1}^{n}|\mathcal{G}_{j}^{n}]),\overline{z}_{j}^{n})\delta +\overline{d}%
_{j}^{n}-\bar{z}_{j}^{n}\varepsilon _{j+1}^{n}\sqrt{\delta },
\label{exp-RBSDE2} \\
\overline{y}_{j}^{n} &\geq &L_{j}^{n},(\overline{y}_{j}^{n}-L_{j}^{n})%
\overline{d}_{j}^{n}=0.  \nonumber
\end{eqnarray}

Apply similar techniques of Lemma \ref{est-expBSDE} to (\ref{exp-RBSDE2}),
we have
\begin{eqnarray*}
E|\bar{y}_{j}^{n}|^{2} &=&E|\overline{y}_{j+1}^{n}|^{2}-E|\bar{z}%
_{j}^{n}|^{2}\delta +2E[\bar{y}_{j+1}^{n}\cdot g(t_{j},E[\overline{y}%
_{j+1}^{n}|\mathcal{G}_{j}^{n}],\overline{z}_{j}^{n})]\delta +2E[\bar{y}%
_{j}^{n}\cdot \overline{d}_{j}^{n}] \\
&&+E|g(t_{j},E[\overline{y}_{j+1}^{n}|\mathcal{G}_{j}^{n}],\overline{z}%
_{j}^{n})|^{2}\delta ^{2}-E[(\overline{d}_{j}^{n})^{2}].
\end{eqnarray*}
In view of $(\overline{y}_{j}^{n}-L_{j}^n)\overline{d}_{j}^{n}=0$ and $%
\overline{d}_{j}^{n}\geq 0$, it follows
\begin{eqnarray*}
E|\bar{y}_{j}^{n}|^{2}+E|\bar{z}_{j}^{n}|^{2}\delta  &\leq &E|\overline{y}%
_{j+1}^{n}|^{2}+2E[\bar{y}_{j+1}^{n}\cdot g(t_{j},E[\overline{y}_{j+1}^{n}|%
\mathcal{G}_{j}^{n}],\overline{z}_{j}^{n})]\delta +2E[(L_{j}^{n})^{+}\cdot
\overline{d}_{j}^{n}] \\
&&+E[|g(t_{j},E[\overline{y}_{j+1}^{n}|\mathcal{G}_{j}^{n}],\overline{z}%
_{j}^{n})|^{2}\delta ^{2}] \\
&\leq &E|\overline{y}_{j+1}^{n}|^{2}+(\delta +3\delta
^{2})E[|g(t_{j},0,0)|^{2}]+(\frac{1}{4}\delta +3\mu ^{2}\delta ^{2})E[(%
\overline{z}_{j}^{n})^{2}] \\
&&+\delta (1+2\mu +4\mu ^{2}+3\mu ^{2}\delta )E|\overline{y}%
_{j+1}^{n}|^{2}+2E[(L_{j}^{n})^{+}\cdot \overline{d}_{j}^{n}]
\end{eqnarray*}
Notice that $3\mu ^{2}\delta <\frac{1}{2}$, since $6\mu ^{2}\delta
<1$. Taking the sum for $j=i,\cdots ,n-1$, it yields
\begin{eqnarray*}
&&E|\bar{y}_{i}^{n}|^{2}+\frac{1}{4}\sum_{j=i}^{n-1}E|\bar{z}%
_{j}^{n}|^{2}\delta  \\
&\leq &E|\xi ^{n}|^{2}+(\delta +3\delta
^{2})E\sum_{j=i}^{n-1}[|g(t_{j},0,0)|^{2}]+\delta (\frac{3}{2}+2\mu +4\mu
^{2})E\sum_{j=i}^{n-1}|\overline{y}_{j+1}^{n}|^{2} \\
&&+\alpha E[\sup_{j}((L_{j}^{n})^{+})^{2}]+\frac{1}{\alpha }%
E[(\sum_{j=i}^{n-1}\overline{d}_{j}^{n})^{2}],
\end{eqnarray*}
where $\alpha $ is a constant to be decided later. Since $\overline{d}%
_{j}^{n}=\overline{y}_{j}^{n}-\overline{y}_{j+1}^{n}-g(t_{j},(E[\overline{y}%
_{j+1}^{n}|\mathcal{G}_{j}^{n}]),\overline{z}_{j}^{n})\delta +\bar{z}%
_{j}^{n}\varepsilon _{j+1}^{n}\sqrt{\delta }$, we get
\[
\sum_{j=i}^{n-1}\overline{d}_{j}^{n}=\overline{y}_{i}^{n}-\xi
^{n}-\sum_{j=i}^{n-1}g(t_{j},(E[\overline{y}_{j+1}^{n}|\mathcal{G}_{j}^{n}]),%
\overline{z}_{j}^{n})\delta +\sum_{j=i}^{n-1}\bar{z}_{j}^{n}\varepsilon
_{j+1}^{n}\sqrt{\delta },
\]
taking square and expectation on both sides, it follows
\begin{eqnarray}
E[(\sum_{j=i}^{n-1}\overline{d}_{j}^{n})^{2}] &\leq &4E|\bar{y}%
_{i}^{n}|^{2}+4E|\xi ^{n}|^{2}+12\delta T
E\sum_{j=i}^{n-1}[|g(t_{j},0,0)|^{2}]+12\mu ^{2}\delta \sum_{j=i}^{n-1}E|%
\overline{y}_{j+1}^{n}|^{2}  \label{est-d} \\
&&+4\delta (3\mu ^{2}+1)\sum_{j=i}^{n-1}|\bar{z}_{j}^{n}|^{2}.
\nonumber
\end{eqnarray}
Set $\alpha =32$, notice that $\delta (3\mu ^{2}+1)<\frac{1}{2}$,
then $\frac{\delta (3\mu ^{2}+1)}{8}<\frac{1}{16}$, we get
\begin{eqnarray*}
\frac{7}{8}E|\bar{y}_{i}^{n}|^{2}&\leq &\frac{11}{8}E|\xi
^{n}|^{2}+(\frac{9}{8}\delta +3\delta
^{2})E\sum_{j=i}^{n-1}[|g(t_{j},0,0)|^{2}]+32E[\sup_{j}((L_{j}^{n})^{+})^{2}]
\\
&&+\delta (\frac{3}{2}+2\mu +\frac{35}{8}\mu ^{2})E\sum_{j=i}^{n-1}|\overline{y}%
_{j+1}^{n}|^{2}.
\end{eqnarray*}
Then apply Lemma \ref{dis-gro}, in view of assumption that implies
$\delta (\frac{3}{2}+2\mu +\frac{35}{8}\mu ^{2})<1$, we obtain
\[
\sup_{j}E[\left| \overline{y}_{j}^{n}\right| ^{2}\leq C_{\xi
^{n},g,L}^{R}.
\]
It follows from the estimations of $\overline{z}_{j}^{n}$ and $%
\overline{d}_{j}^{n}$ that
\[
E[\sum_{j=0}^{n-1}\left|
\overline{z}_{j}^{n}\right| ^{2}\delta +(\sum_{j=0}^{n-1}\overline{d}%
_{j}^{n})^{2}]\leq C_{\xi ^{n},g,L}^{R}
\]
As Lemma \ref{est-expBSDE}, using Burkholder-Davis-Gundy inequality
and similar techniques, we get the results. $%
\square $

Then we have following convergence result for explicit reflected scheme.

\begin{theorem}
Under the same assumptions of Theorem \ref{conv-imp}, the discrete solutions
$\{(\bar{Y}^{n},\bar{Z}^{n})\}_{n=1}^{\infty }$ of the explicit reflected
scheme converges to the solution $(Y,Z)$ of (\ref{RBSDE1b1}) in the
following senses: as $n\rightarrow \infty $%
\begin{equation}
E[\sup_{0\leq t\leq T}\left| \bar{Y}_{t}^{n}-Y_{t}\right|
^{2}+\int_{0}^{T}\left| \bar{Z}_{s}^{n}-Z_{s}\right|
^{2}ds]\rightarrow 0. \label{conv-exp-R}
\end{equation}
Moreover $E[\sup_{0\leq t\leq T}(K_{t}-\overline{K}_{t}^{n})^{2}]\rightarrow 0$%
.
\end{theorem}

\noindent\textbf{Proof. }Thanks to convergence results of Theorem
\ref{conv-imp}, it suffices to prove
\begin{equation}
E[\sup_{0\leq t\leq T}\left| Y_{t}^{n}-\bar{Y}_{t}^{n}\right|
^{2}+\int_{0}^{T}\left| Z_{s}^{n}-\bar{Z}_{s}^{n}\right|
^{2}ds]\rightarrow 0. \label{conv-exp-R1}
\end{equation}
Recall the implicit reflected scheme and explicit reflected scheme: for $%
0\leq j\leq n-1$,
\begin{eqnarray*}
y_{j}^{n} &=&y_{j+1}^{n}+g(t_{j},y_{j}^{n},z_{j}^{n})\delta
+d_{j}^{n}-z_{j}^{n}\varepsilon _{j+1}^{n}\sqrt{\delta }, \\
\overline{y}_{j}^{n} &=&\overline{y}_{j+1}^{n}+g(t_{j},E[\overline{y}%
_{j+1}^{n}|\mathcal{G}_{j}^{n}]),z_{j}^{n})\delta +\overline{d}_{j}^{n}-%
\overline{z}_{j}^{n}\varepsilon _{j+1}^{n}\sqrt{\delta },
\end{eqnarray*}
Consider the difference, we have
\begin{eqnarray*}
&&E\left| y_{j}^{n}-\overline{y}_{j}^{n}\right| ^{2} \\
&=&E\left| y_{j+1}^{n}-\overline{y}_{j+1}^{n}\right| ^{2}-\delta E\left|
z_{j}^{n}-\overline{z}_{j}^{n}\right| ^{2}+2\delta E[(y_{j}^{n}-\overline{y}%
_{j}^{n})(g(t_{j},y_{j}^{n},z_{j}^{n})-g(t_{j},E[\overline{y}_{j+1}^{n}|%
\mathcal{G}_{j}^{n}],\overline{z}_{j}^{n}))] \\
&&+2E[(y_{j}^{n}-\overline{y}_{j}^{n})(d_{j}^{n}-\overline{d}%
_{j}^{n})]-\delta ^{2}E[(g(t_{j},y_{j}^{n},z_{j}^{n})-g(t_{j},E[\overline{y}%
_{j+1}^{n}|\mathcal{G}_{j}^{n}],\overline{z}_{j}^{n}))^{2}] \\
&&-2\delta E[(d_{j}^{n}-\overline{d}%
_{j}^{n})(g(t_{j},y_{j}^{n},z_{j}^{n})-g(t_{j},E[\overline{y}_{j+1}^{n}|%
\mathcal{G}_{j}^{n}],\overline{z}_{j}^{n}))]-E\left| d_{j}^{n}-\overline{d}%
_{j}^{n}\right| ^{2} \\
&\leq &E\left| y_{j+1}^{n}-\overline{y}_{j+1}^{n}\right| ^{2}-\delta E\left|
z_{j}^{n}-\overline{z}_{j}^{n}\right| ^{2}+2\delta E[(y_{j}^{n}-\overline{y}%
_{j}^{n})(g(t_{j},y_{j}^{n},z_{j}^{n})-g(t_{j},E[\overline{y}_{j+1}^{n}|%
\mathcal{G}_{j}^{n}],\overline{z}_{j}^{n}))],
\end{eqnarray*}
in view of $-a^{2}-2ab-b^{2}=-(a+b)^{2}\leq 0$ and
\begin{eqnarray*}
(y_{j}^{n}-\overline{y}_{j}^{n})(d_{j}^{n}-\overline{d}_{j}^{n})
&=&(y_{j}^{n}-L_{j}^{n})d_{j}^{n}+(\overline{y}_{j}^{n}-L_{j}^{n})(\overline{%
d}_{j}^{n}) \\
&&-(\overline{y}_{j}^{n}-L_{j}^{n})d_{j}^{n}-(\overline{y}%
_{j}^{n}-L_{j}^{n})(d_{j}^{n}) \\
&\leq &0.
\end{eqnarray*}
We take sum over $j$ from $i$ to $n-1$, with $\xi ^{n}-\overline{\xi }^{n}=0$%
, then get
\[
E\left| y_{i}^{n}-\overline{y}_{i}^{n}\right| ^{2}+\delta
\sum_{j=i}^{n-1}E\left| z_{j}^{n}-\overline{z}_{j}^{n}\right| ^{2}\leq
2\delta \sum_{j=i}^{n-1}E[(y_{j}^{n}-\overline{y}%
_{j}^{n})(g(t_{j},y_{j}^{n},z_{j}^{n})-g(t_{j},E[\overline{y}_{j+1}^{n}|%
\mathcal{G}_{j}^{n}],\overline{z}_{j}^{n}))].
\]
Now we are in the same situation as in the proof of Theorem \ref{conv1}. By
similar methods, with Lemma \ref{est-expRBSDE}, and
\begin{eqnarray*}
&&2\mu \delta E[\left| y_{j}^{n}-\overline{y}_{j}^{n}\right| \cdot \left|
y_{j}^{n}-E[\overline{y}_{j+1}^{n}|\mathcal{G}_{j}^{n}]\right| ] \\
&=&2\mu \delta E[\left| y_{j}^{n}-\overline{y}_{j}^{n}\right| \cdot \left|
y_{j}^{n}-\overline{y}_{j}^{n}+g(t_{j},E[\overline{y}_{j+1}^{n}|\mathcal{G}%
_{j}^{n}]),z_{j}^{n})\delta +\overline{d}_{j}^{n}\right| ] \\
&\leq &(2\mu +1)\delta E[\left| y_{j}^{n}-\overline{y}_{j}^{n}\right|
^{2}]+2\mu ^{2}\delta E[\delta ^{2}\left| g(t_{j},E[\overline{y}_{j+1}^{n}|%
\mathcal{G}_{j}^{n}]),z_{j}^{n})\right| ^{2}+(\overline{d}_{j}^{n})^{2}],
\end{eqnarray*}
we obtain
\begin{equation}
E\left| y_{i}^{n}-\overline{y}_{i}^{n}\right| ^{2}+\frac{\delta }{2}%
\sum_{j=i}^{n-1}E\left| z_{j}^{n}-\overline{z}_{j}^{n}\right| ^{2}\leq (2\mu
^{2}+2\mu +1)\delta \sum_{j=i}^{n-1}E\left| y_{j}^{n}-\overline{y}%
_{j}^{n}\right| ^{2}+\delta C_{\xi ^{n},g,L}^{R},  \label{est-exp-d}
\end{equation}
where $C_{\xi ^{n},g,L}^{R}$ is a constant only depends on $\xi ^{n}$, $%
g(\cdot ,0,0)$, $\mu $ and $L$. By Lemma \ref{dis-gro}, we get
\[
\sup_{j\leq n}E\left| y_{j}^{n}-\overline{y}_{j}^{n}\right| ^{2}\leq C\delta
e^{(2\mu +2\mu ^{2}+1)T}.
\]
From (\ref{est-exp-d}), it follows (\ref{conv-exp-R1}), which
implies $\lim_{n\rightarrow \infty} \delta  %
\sum_{j=i}^{n-1}E\left| z_{j}^{n}-\overline{z}_{j}^{n}\right| ^{2}=
0$. Then (\ref {conv-exp-R}) follows by  using
Burkholder-Davis-Gundy inequality, similar techniques and
estimations results from Lemma \ref{est-expRBSDE}. In fact, we get
\[
E[\sup_{j\leq n}\left| y_{j}^{n}-\overline{y}_{j}^{n}\right|
^{2}]\leq C_{\mu}  \delta  %
\sum_{j=i}^{n-1}E\left| y_{j}^{n}-\overline{y}_{j}^{n}\right| ^{2}+\delta  %
\sum_{j=i}^{n-1}E\left| z_{j}^{n}-\overline{z}_{j}^{n}\right| ^{2}.
\]

For the convergence of $\overline{K}^{n}$, for $0\leq t\leq T$, since
\begin{eqnarray*}
K_{t}^{n} &=&Y_{0}^{n}-Y_{t}^{n}-\sum_{j=0}^{[t/\delta
]}g(t_{j},y_{j}^{n},z_{j}^{n})\delta +\sum_{j=0}^{[t/\delta
]}z_{j}^{n}\varepsilon _{j+1}^{n}\sqrt{\delta }, \\
\overline{K}_{t}^{n} &=&\overline{Y}_{0}^{n}-\overline{Y}_{t}^{n}-%
\sum_{j=0}^{[t/\delta ]}g(t_{j},E[\overline{y}_{j+1}^{n}|\mathcal{G}%
_{j}^{n}],\overline{z}_{j}^{n})\delta +\sum_{j=0}^{[t/\delta ]}\overline{z}%
_{j}^{n}\varepsilon _{j+1}^{n}\sqrt{\delta },
\end{eqnarray*}
with Lipschitz condition of $g$ and BDG inequality, we get
\begin{eqnarray*}
E[\sup_{0\leq t\leq T}(K_{t}^{n}-\overline{K}_{t}^{n})^{2}] &\leq &4\left| Y_{0}^{n}-\overline{Y}%
_{0}^{n}\right| ^{2}+4E[\sup_{0\leq t\leq T}\left|
Y_{t}^{n}-\overline{Y}_{t}^{n}\right|
^{2}]+8\delta \mu ^{2}\sum_{j=0}^{[t/\delta ]}\left| y_{j}^{n}-E[\overline{y}%
_{j+1}^{n}|\mathcal{G}_{j}^{n}]\right| ^{2} \\
&&+4\delta (2\mu ^{2}+c_2)\sum_{j=0}^{[t/\delta ]}\left| z_{j}^{n}-\overline{z}%
_{j}^{n}\right| ^{2} \\
&=&4\left| Y_{0}^{n}-\overline{Y}_{0}^{n}\right| ^{2}+4E[\sup_{0\leq t\leq T}\left| Y_{t}^{n}-%
\overline{Y}_{t}^{n}\right| ^{2}]+4\delta (2\mu
^{2}+c_2)\sum_{j=0}^{[t/\delta
]}\left| z_{j}^{n}-\overline{z}_{j}^{n}\right| ^{2} \\
&&+24\delta \mu ^{2}\sum_{j=0}^{[t/\delta ]}\left| y_{j}^{n}-\overline{y}%
_{j}^{n}\right| ^{2}+\delta C_{\xi ^{n},g,L}^{R}.
\end{eqnarray*}
From (\ref{conv-exp-R1}) and convergence of $K_{t}^{n}$ to $K_{t}$ in $%
\mathbf{L}^{2}(\mathcal{F}_{t})$, we obtain the convergence of $\overline{K}%
_{t}^{n}$ to $K_{t}$. $\square $

\subsection{Simulations of Reflected BSDEs with one lower barrier}

For calculation convenience, we consider the case when $T=1$, and
begin with $y_{n}^{n}=\xi ^{n}$ backwardly solve
$(y_{j}^{n},z_{j}^{n},d_{j}^{n})$, for $j=n-1,\cdots 1,0$. Consider
the amount of total calculation for most general case, we only treat
a very simple situation: $\xi =\phi (B_{1})$, $L_{t}=\psi (t,B(t))$,
where $\phi $
and $\psi $ are real regular functions defined on $\mathbb{R }$ and $%
[0,1]\times \mathbb{R}$ respectively. As for BSDE, we have also
developed a Matlab toolbox for calculating and simulating solutions
of reflected BSDEs. This toolbox is similar to the one for BSDE and
can be downloaded from
http://159.226.47.50:8080/iam/xumingyu/English.jsp, by clicking
'Preprint' on the left side.

Here we consider following case: $g(t,y,z)=-\left| y+z\right| $, $\xi =\Phi
(B_{1})=2\sin (B_{1})$, $L_{t}=\Psi (t,B_{t})=\sin (B_{t}+\frac{\pi }{2})-2$
and $n=400$.

After inputting the parameters, we run the calculation program using
reflected explicit scheme, then get all prossible results of $y$. We may
notice that at $t=1$, $\xi \geq L_{1}$ does not always hold. But the
numerical scheme still works as well. In fact, in such case the increasing
process $K$ as well as $y$ has a jump of size $(L_{1}-\xi )^{+}$ at $t=1$,
which pushes the solution $y_{t-}$, i.e. $y_{n-1}$ in our case, to stay
above the barrier $L$. Then both $K$ and $y$ act as the terminal condition
is $(\xi -L_{1})^{+}+L_{1}$, which is always bigger than $L_{1}$.

Now we will see some properties of the trajectory of solution $y$ in
the Figure 3. In the upper portion of Figure 3, the below surface
shows the barrier $L$ in 3-dimensional, as well the upper one is for
the solution $y$. Then we use programs to generate two trajectories
of the discrete Brownian motion $(B_{j}^{n,i})_{0\leq j\leq n}$, for
$i=1,2$, which are drawn on horizontal plane. The value of
$y_{j}^{n,i}(i=1,2)$ with respect to these Brownian samples, are
showed on the solution surface, and we use the fine vertical line to
give correspondence between two group of trajectories of $y$ and
$B$. The remainder of the figure shows respectively the trajectories
of the force $K_{j}^{n,i}=\sum_{k=0}^{j}d_{k}^{n,i}\;(i=1,2)$
corresponding to the value of $y_{j}^{n,i}(i=1,2)$, and
$y_{j}^{n,i}-L_{j}^{n,i}(i=1,2)$.

In the upper portion we can see that there is an area where two
surfaces (the solution surface and the barrier surface) stick
together. When the
trajectory of solution $y_{j}^{n}$ goes into this area, the force $%
K_{j}^{n,i}$ will push $y_{j}^{n}$ upward. Indeed, if we don't have
the barrier here, $y_{j}^{n,i}$ intends becoming smaller than the
reflecting
barrier $L_{j}^{n,i}$, so to keep $y_{j}^{n,i}$ being no less than $%
L_{j}^{n,i}$, the action of forces $K_{j}^{n,i}$ are necessary.
Comparing these two trajectories, we can see that one trajectory,
noted as $K_{j}^{n,1} $ pushes upwards the corresponding trajectory
of solution $y_{j}^{n,1}$, while the other one noted as
$K_{j}^{n,2}$, keeps zero, since $y_{j}^{n,1}$ goes into the
sticking area but the trajectory $y_{j}^{n,2}$ with respect to
$K_{j}^{n,2}$ does not.
\begin{center}
\centering
\includegraphics[height=125mm,width=155mm]{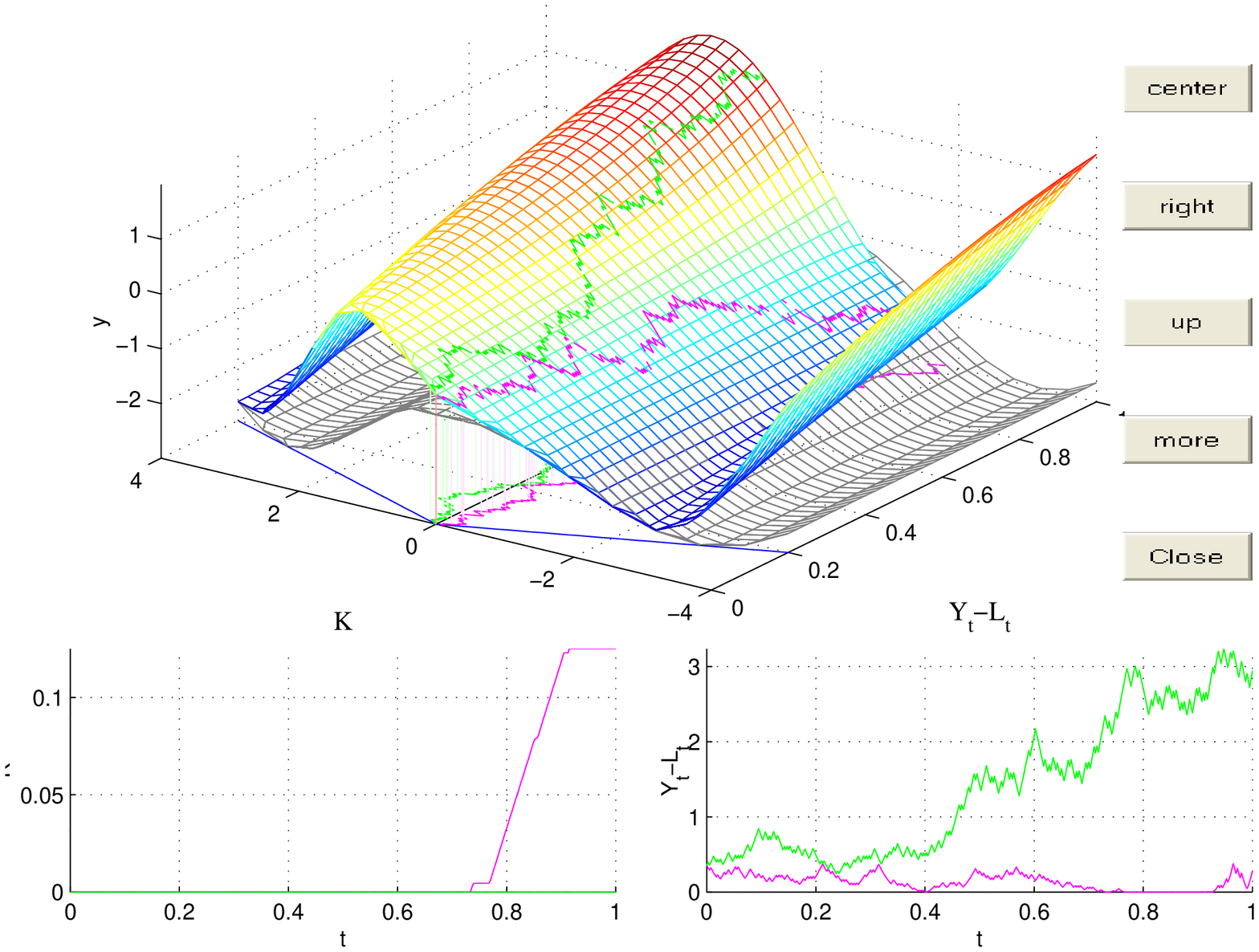}\\
Figure 3: The solution on surface
\end{center}

Compare the two sub-figures, which are below the main sub-figure, we can
easily find out that the $K_{j}^{n,i}$ only increases when $%
y_{j}^{n,i}-L_{j}^{n,i}$ takes the value $0$; but the converse is not always
true, when $y_{j}^{n,i}-L_{j}^{n,i}=0$, $K_{j}^{n,i}$ does not necessary
increase.

About this point, we can also see in Figure 4.This figure shows two
groups
of 3-dimensional dynamic trajectories $%
(t_{j},B_{j}^{n,i},Y_{j}^{n,i})_{(i=1,2)}$ and $%
(t_{j},B_{j}^{n,i},Z_{j}^{n,i})_{(i=1,2)}$ and, simultaneously, two
groups
of 2-dimensional trajectories of $(t_{j},Y_{j}^{n,i})_{(i=1,2)}$ and $%
(t_{j},Z_{j}^{n,i})_{(i=1,2)}$. For remainder sub-figures, the
above-right one is for the trajectories $K_{j}^{n,i}(i=1,2)$, and
while the below-left one is for $y_{j}^{n,i}-L_{j}^{n,i}(i=1,2)$,
then comparing the these two
sub-figures, as in Figure 3, we can see clearly the relation between $%
K_{j}^{n,i}(i=1,2)$ and $y_{j}^{n,i}-L_{j}^{n,i}(i=1,2)$. Moreover
one trajectory, noted as $K^{n,1}$ as well as $y^{n,1}$jumps at
$t=1$, since its terminal value is less than the barrier.

\begin{center}
\centering
\includegraphics[height=110mm,width=150mm]{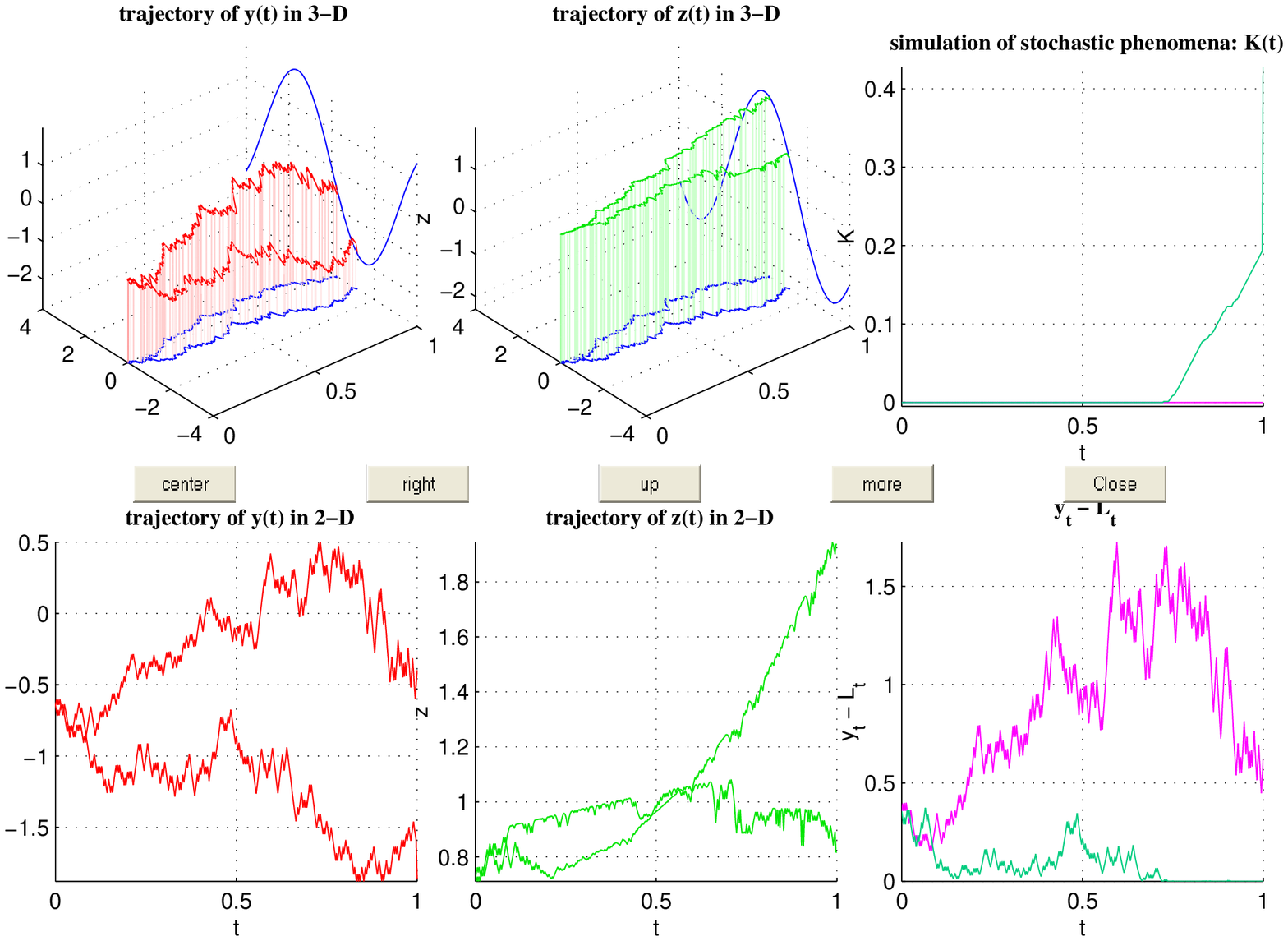}\\
Figure 4: Simulation of $y$, $z$, $K$
\end{center}

Then we list out some numerical results for the reflected scheme and
explicit implicit penalization scheme, and we can see that as the penalized
parameter $p$ converge to infinity, $y_{0}^{p,n}$ converge to $y_{0}^{n}$.
Consider the same parameters as above: $f(y,z)=-\left| y+z\right| $, $\xi
=\Phi (B_{1})=2\sin (B_{1})$, $L_{t}=\Psi (t,B_{t})=\sin (B_{t}+\frac{\pi }{2%
})-2$. Then as the following tablet showing:
\[
\begin{array}{lll}
n=400, & \mbox{reflected explicit scheme: } & y_{0}^{n}=-0.6430, \\
& \mbox{penalization scheme:} &
\begin{tabular}{|c|c|c|c|c|}
\hline
$p$ & $20$ & $200$ & $2000$ & $2\times 10^{4}$ \\ \hline
$y_{0}^{p,n}$ & $-0.6553$ & $-0.6444$ & $-0.6431$ & $-0.6430$ \\ \hline
\end{tabular}
\end{array}
\]
\[
\begin{array}{lll}
n=1000, & \mbox{reflected explicit scheme:} & y_{0}^{n}=-0.6425, \\
& \mbox{penalization scheme:} &
\begin{tabular}{|c|c|c|c|c|}
\hline
$p$ & $20$ & $200$ & $2000$ & $2\times 10^{4}$ \\ \hline
$y_{0}^{p,n}$ & $-0.6550$ & $-0.6441$ & $-0.6427$ & $-0.6425$ \\ \hline
\end{tabular}
\end{array}
\]
\[
\begin{array}{lll}
n=2000, & \mbox{reflected explicit scheme:} & y_{0}^{n}=-0.6424, \\
& \mbox{penalization scheme:} &
\begin{tabular}{|c|c|c|c|c|}
\hline
$p$ & $20$ & $200$ & $2000$ & $2\times 10^{4}$ \\ \hline
$y_{0}^{p,n}$ & $-0.6549$ & $-0.6439$ & $-0.6426$ & $-0.6424$ \\ \hline
\end{tabular}
\end{array}
\]

\begin{remark}
For BSDE with two reflecting barriers, we introduced also reflected implicit
and explicit scheme as well as penalization schemes. The proofs of
convergence and simulations results can be found in \cite{X2007}.
\end{remark}

\section{$\Gamma $-constrained BSDEs}

In this section, we consider a $g$-supersolution $(Y,Z,A)$ with constraint $%
(Y_{t},Z_{t})\in \Gamma _{t}$ of the following form:
\begin{equation}
Y_{t}=\xi +\int_{t}^{T}g(s,Y_{s},Z_{s})ds+A_{T}-A_{t}-\int_{t}^{T}Z_{s}dB_{s}
\label{CBSDE}
\end{equation}
with $d_{\Gamma }(Y_{t},Z_{t})=0,$ a.e., a.s., where $\Gamma $ is a
nonempty closed subset of $\mathbb{R\times R}$ and $d_{\Gamma }$ is
the distance function of $\Gamma $, i.e., $d_{\Gamma
}(y,z)=\inf_{(y^{\prime },z^{\prime })\in \Gamma }\{|y-y^{\prime
}|+|z-z^{\prime }|\}$. It is clear that $\Gamma $ is a Lipschitz
function
\[
|d_{\Gamma }(y,z)-d_{\Gamma }(y^{\prime },z^{\prime })|\leq |y-y^{\prime
}|+|z-z^{\prime }|.
\]
Such a $g$-supersolution $(Y,Z,A)$ is called a $\Gamma $-constrained $g$%
-supersolution.

As before we assume that $g$ satisfies Lipschitz condition (\ref{Lip}) and
that assumption \ref{term-num} holds for $\xi $. From \cite{Peng99}, we have
the existence of smallest solutions for (\ref{CBSDE}):

\begin{theorem}
\label{constr}If there exists at least one $\Gamma $-constrained $g$%
-supersolution of (\ref{CBSDE}), then the equation admits a smallest $\Gamma
$-constrained $g$-supersolution $(Y,Z,A)$. Moreover, $(Y,Z)$ is the limit of
the following sequence of penalization solutions:
\begin{equation}
Y_{t}^{p}=\xi +\int_{t}^{T}g(s,Y_{s}^{p},Z_{s}^{p})ds+p\int_{t}^{T}d_{\Gamma
}(Y_{s}^{p},Z_{s}^{p})ds-\int_{t}^{T}Z_{s}^{p}dB_{s},  \label{pCBSDE}
\end{equation}
in the sense of
\[
\lim_{p\rightarrow \infty
}E\int_{0}^{T}[|Y_{t}-Y_{t}^{p}|^{2}+|Z_{t}-Z_{t}^{p}|^{\beta }]dt=0,\ 1\leq
\beta <2.
\]
\end{theorem}

This smallest $\Gamma $-constrained $g$-supersolution is called a $g_{\Gamma
}$-solution. And such equation can be considered as a BSDE with singular
coefficient $g_{\Gamma }:=g(t,y,z)+\infty \cdot d_{\Gamma }(y,z)$. It easy
to check that when $\xi ^{+}\in L^{\infty }(\mathcal{F}_{T})$, and there
exists a large enough constant $C_{0}$ such that for $y\geq C_{0}$%
\[
g(t,y,0)\leq C_{0}+\mu |y|,\;\ \text{and \ }(y,0)\in \Gamma ,
\]
there exists a $\Gamma $-constrained $g$-supersolution of (\ref{CBSDE})(see
Peng and Xu \cite{PX2005}). Then by Theorem \ref{constr}, a $g_{\Gamma }$%
-solution exists. In this section we will work under these assumptions. We
now derive a numerical scheme applying convergence results in Theorem \ref
{constr}.

\subsection{Constraint on $Z$}

First we consider the case when constraint is only on process $Z$
and invariant in $t$ , i.e. $\Gamma $ is a close subset in
$\mathbf{R}$. And we require $Z\in \Gamma $, i.e. $d_{\Gamma
}(Z_{t})=0$, a.e.a.s.. After same discretization for BSDEs
introduced in Section 2, for each positive number $p
$ we have the following penalization discrete equation on small interval $%
[j\delta ,(j+1)\delta ]$%
\begin{equation}
y_{j}^{p,n}=y_{j+1}^{p,n}+g(t_{j},y_{j}^{p,n},z_{j}^{p,n})\delta +pd_{\Gamma
}(z_{j}^{p,n})\delta -z_{j}^{p,n}\sqrt{\delta }\varepsilon _{j+1},
\label{dis-pc}
\end{equation}
with discrete terminal condition: $y_{n}^{n}:=\xi ^{n}$.

Now we need to find a way to find $\mathcal{G}_{j}^{n}$-measurable $%
(y_{j}^{p,n},z_{j}^{n,p})$ to satisfy (\ref{dis-pc}) with $%
(y_{j+1}^{p,n},z_{j+1}^{p,n})$. It is easy to get $%
z_{j}^{p,n}=\frac{1}{\sqrt{\delta}}E[y_{j+1}^{p,n}\varepsilon _{j+1}^{n}|\mathcal{F}_{j}^{n}]=\frac{%
1}{2\sqrt{\delta }}(y_{j+1}^{p,n}|_{\varepsilon
_{j}^{n}=1}-y_{j+1}^{p,n}|_{\varepsilon _{j}^{n}=-1})$. Substitute it into (%
\ref{dis-pc}), it follows a equation of $y_{j}^{p,n}$as
\[
y_{j}^{p,n}=E[y_{j+1}^{p,n}|\mathcal{G}%
_{j}^{n}]+g(t_{j},y_{j}^{p,n},z_{j}^{p,n})\delta +pd_{\Gamma
}(z_{j}^{p,n})\delta .
\]
So apply the implicit scheme for BSDE in Section 2, we get
\[
y_{j}^{p,n}=\Theta ^{-1}(E[y_{j+1}^{p,n}|\mathcal{G}_{j}^{n}]+pd_{\Gamma
}(z_{j}^{p,n})\delta ),
\]
where $\Theta (y)=y-g(t_{j},y,z_{j}^{p,n})\delta $. While the explicit
method gives
\[
\overline{y}_{j}^{p,n}=E[\overline{y}_{j+1}^{p,n}|\mathcal{G}%
_{j}^{n}]+g(t_{j},E[\overline{y}_{j+1}^{p,n}|\mathcal{G}_{j}^{n}],\overline{z%
}_{j}^{p,n})\delta +pd_{\Gamma }(\overline{z}_{j}^{p,n})\delta .
\]
The interesting point here is that the penalization of $z^{p,n}$ with
respect to $z$ is not directly on $z^{p,n}$, it act on $y^{n,p}$ to
influence $z^{p,n}$.

We have

\begin{theorem}[Convergence Theorem]
Define
\[
Y_{t}^{p,n}=y_{[t/\delta ]}^{p,n},Z_{t}^{p,n}=z_{[t/\delta ]}^{p,n},%
\overline{Y}_{t}^{p,n}=\overline{y}_{[t/\delta ]}^{p,n},\overline{Z}%
_{t}^{p,n}=\overline{z}_{[t/\delta ]}^{p,n}.
\]
Here $y_{j}^{p,n}$ and $\overline{y}_{j}^{p,n}$, $0\leq j\leq n$,
can come from either implicit scheme or explicit scheme
respectively. Under assumption \ref{term-num}, and $g$ satisfying
Lipschitz condition. Then
\[
\lim_{p\rightarrow\infty}\lim_{n\rightarrow
\infty}E[\int_{0}^{T}\left| Y_{s}^{p,n}-Y_{s}\right|
^{2}ds+\int_{0}^{T}\left| Z_{s}^{p,n}-Z_{s}\right| ^{\beta }ds]=
0,1\leq \beta <2.
\]
\end{theorem}

\noindent \textbf{Proof. }By Theorem \ref{constr}, for any $\varepsilon >0$,
there exists $p_{0}>0$ such that for each $p>p_{0}$,
\[
E[\int_{0}^{T}\left| Y_{s}^{p}-Y_{s}\right| ^{2}ds^{2}+\int_{0}^{T}\left|
Z_{s}^{p}-Z_{s}\right| ^{\beta }ds]\leq \varepsilon .
\]
Moreover, by Theorem \ref{conv1}, for implicit scheme, we have as $%
n\rightarrow \infty $%
\[
E[\sup_{0\leq t\leq T}\left| Y_{t}^{p_{0},n}-Y_{t}^{p_{0}}\right|
^{2}+\int_{0}^{T}\left| Z_{s}^{p_{0},n}-Z_{s}^{p_{0}}\right|
^{2}ds]\rightarrow 0.
\]
For explicit scheme, the result follows from
\[
\sup_{0\leq t\leq T}E[\left| \overline{Y}_{t}^{p_{0},n}-Y_{t}^{p_{0}}\right|
^{2}+\int_{0}^{T}\left| \overline{Z}_{s}^{p_{0},n}-Z_{s}^{p_{0}}\right|
^{2}ds]\rightarrow 0.
\]
\noindent $\square $

To illustrate calculation and simulation in our software package, we
consider the case $\Gamma =[a,b]$ with $a\leq 0\leq b$. Then $d_{\Gamma
}(z)=(z-a)^{-}+(z-b)^{+}$. The default setting is $g(t,y,z)=-2\left|
y+z\right| -1$, $\xi =\left| B_{1}\right| $, with $a=-0.5$, $b=0.8$, $p=20$
and $n=400$. The surface $u=u^{p,n}(t,x)$ and $v=v^{n,p}(t,x)$ are given
with dynamic simulation $Y_{t}^{p,n}=u^{p,n}(t,B_{t}^{n})$ and $%
Z_{t}^{p,n}=v^{p,n}(t,B_{t}^{n})$ as shown as the upper part of
Figure 5. The lower part of the figure displays the simultaneous
action of the process $A_{t}^{p,n}=p\sum_{j\leq [t/\delta
]}d_{\Gamma }(z_{j}^{p,n})$. The effect of increases in the
$A_{t}^{p,n}$ when $Z^{p,n}$ is less than $-0.5$ and larger than
$0.8$ are clearly shown in Figure 6. But it seems that the
solution is still too sensitive to the choice of $p$ and $n$. If $p\sqrt{%
\delta }>1$, then the numerical solution will explode.

\begin{center}
\centering\includegraphics[totalheight=115mm]{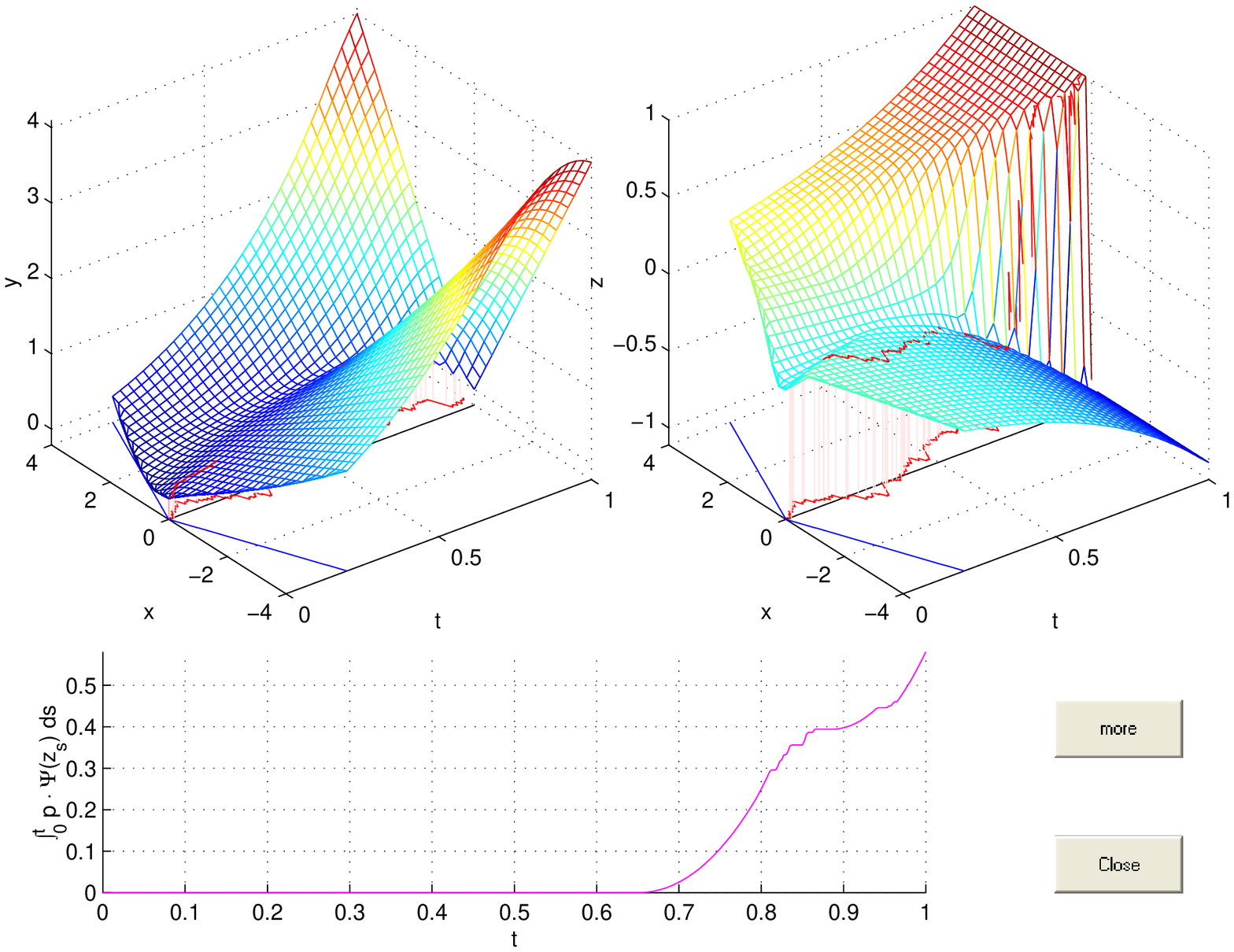} \\[0pt]
Figure 5: The solution surface of BSDE (\ref{dis-pc})

\centering\includegraphics[height=80mm,width=160mm]{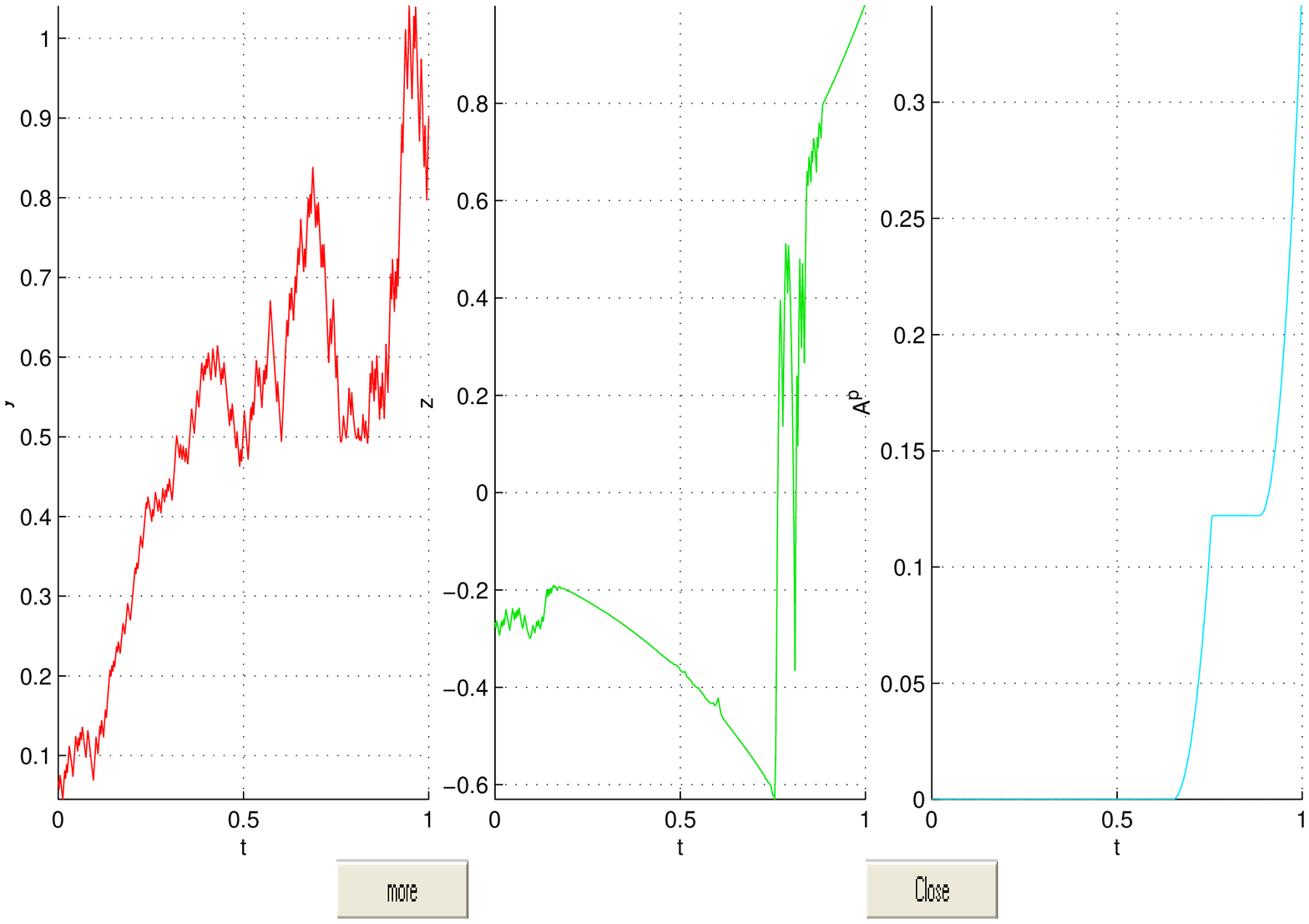} \\[0pt]
Figure 6: A trajectory of solutions of (\ref{dis-pc})
\end{center}

\subsection{BSDE reflected on process $Z$}

Now we consider another special case, when the constraint is $d_{\Gamma
}(y,z)=(y-\phi (z))^{-}$, in other words, we require $y-\phi (z)\geq 0$.
After the same discretization of the time interval, we have following
discrete penalization equation for some $p$ large enough, on the small
interval $[j\delta ,(j+1)\delta ]$, $0\leq j\leq n-1$%
\[
y_{j}^{p,n}=y_{j+1}^{p,n}+g(t_{j},y_{j}^{p,n},z_{j}^{p,n})\delta +p(\phi
(z_{j}^{p.n})-y_{j}^{p,n})^{+}\delta -z_{j}^{p,n}\delta \varepsilon _{j+1}.
\]
Similarly, $z_{j}^{p,n}=\frac{1}{\sqrt{\delta }}E[y_{j+1}^{p,n}\varepsilon _{j+1}^{n}|\mathcal{F}%
_{j}^{n}]=\frac{1}{2\sqrt{\delta }}(y_{j+1}^{p,n}|_{\varepsilon
_{j}^{n}=1}-y_{j+1}^{p,n}|_{\varepsilon _{j}^{n}=-1})$. Then
$y_{j}^{p,n}$ satisfies
\[
y_{j}^{p,n}=E[y_{j+1}^{p,n}|\mathcal{G}%
_{j}^{n}]+g(t_{j},y_{j}^{p,n},z_{j}^{p,n})\delta +p(\phi
(z_{j}^{p.n})-y_{j}^{p,n})^{+}\delta .
\]
Set $\Theta (y)=y-(g(t_{j},y,z_{j}^{p,n})\delta +p(\phi
(z_{j}^{p,n})-y)^{+}\delta )$, with $E[y_{j+1}^{p,n}|\mathcal{G}_{j}^{n}]=%
\frac{1}{2}(y_{j+1}^{p,n}|_{\varepsilon
_{j}^{n}=1}+y_{j+1}^{p,n}|_{\varepsilon _{j}^{n}=-1})$, then our implicit
scheme is given by solving following equation
\[
y_{j}^{p,n}=\Theta ^{-1}(E[y_{j+1}^{p,n}|\mathcal{G}_{j}^{n}]).
\]
Meanwhile, we have also explicit-implicit scheme, which is
\begin{eqnarray*}
\overline{y}_{j}^{p,n} &=&E[\overline{y}_{j+1}^{p,n}|\mathcal{G}%
_{i}^{n}]+g(t_{j},E[\overline{y}_{j+1}^{p,n}|\mathcal{G}_{j}^{n}],\overline{z%
}_{j}^{p,n})\delta  \\
&&+\frac{p\delta }{1+p\delta }(E[\overline{y}_{j+1}^{p,n}|\mathcal{G}%
_{j}^{n}]+g(t_{j},E[\overline{y}_{j+1}^{p,n}|\mathcal{G}_{j}^{n}],\overline{z%
}_{j}^{p,n})\delta -\phi (\overline{z}_{j}^{p,n}))^{-}
\end{eqnarray*}

As previous subsection, we have convergence results of these two schemes.

\begin{theorem}
Define $Y_{t}^{p,n}=y_{[t/\delta ]}^{p,n}$, $Z_{t}^{p,n}=z_{[t/\delta
]}^{p,n}$ and $\overline{Y}_{t}^{p,n}=\overline{y}_{[t/\delta ]}^{p,n}$, $%
\overline{Z}_{t}^{p,n}=\overline{z}_{[t/\delta ]}^{p,n}$. Then we
have, for $1\leq \beta <2$,
\begin{eqnarray*}
\lim_{p\rightarrow
\infty}\lim_{n\rightarrow\infty}E[\int_{0}^{T}\left|
Y_{s}^{p,n}-Y_{s}\right| ^{2}ds+\int_{0}^{T}\left|
Z_{s}^{p,n}-Z_{s}\right| ^{\beta}ds] &= &0, \\
\lim_{p\rightarrow
\infty}\lim_{n\rightarrow\infty}E[\int_{0}^{T}\left|
\overline{Y}_{s}^{p,n}-Y_{s}\right| ^{2}ds+\int_{0}^{T}\left|
\overline{Z}_{s}^{p,n}-Z_{s}\right| ^{\beta}ds] &= &0.
\end{eqnarray*}
\end{theorem}

\noindent\textbf{Proof.} The results follow from Theorem
\ref{constr} and Proposition \ref{ThCovIE}, so we omit the proof.
$\square $

Now we do simulations by explicit-implicit scheme. We consider the case $%
g=-2\left| y+z\right| -1$, $\xi =\left| B_{1}\right| $, $\phi (z)=1.25\times
z$, with penalization parameter $p=10$, and discretization number $n=400$.
In figure 7, we see the surface of solution $Y^{p,n}$ with a trajectory of $%
Y^{p,n}$ on the surface in upper portion, while in two lower subfigures
there presents the trajectory $A^{p,n}=p\sum_{j\leq [t/\delta
]}(y_{j}^{p,n}-\phi (z_{j}^{p,n}))^{-}\delta $ and $y_{j}^{p,n}-\phi
(z_{j}^{p,n})$ corresponding to the one on the surface. We can see that when
$y_{j}^{p,n}-\phi (z_{j}^{p,n})$ is positive, the penalization term will not
work to the process $y_{j}^{p,n}$. About this point we can see more clear in
Figure 8, which presents trajectories of $y_{j}^{p,n}$, $z_{j}^{p,n}$, $%
A^{p,n}$ and $y_{j}^{p,n}-\phi (z_{j}^{p,n})$ in 3 or 2-dimensional
subfigures.

\begin{center}
\centering\includegraphics[height=100mm,width=160mm]{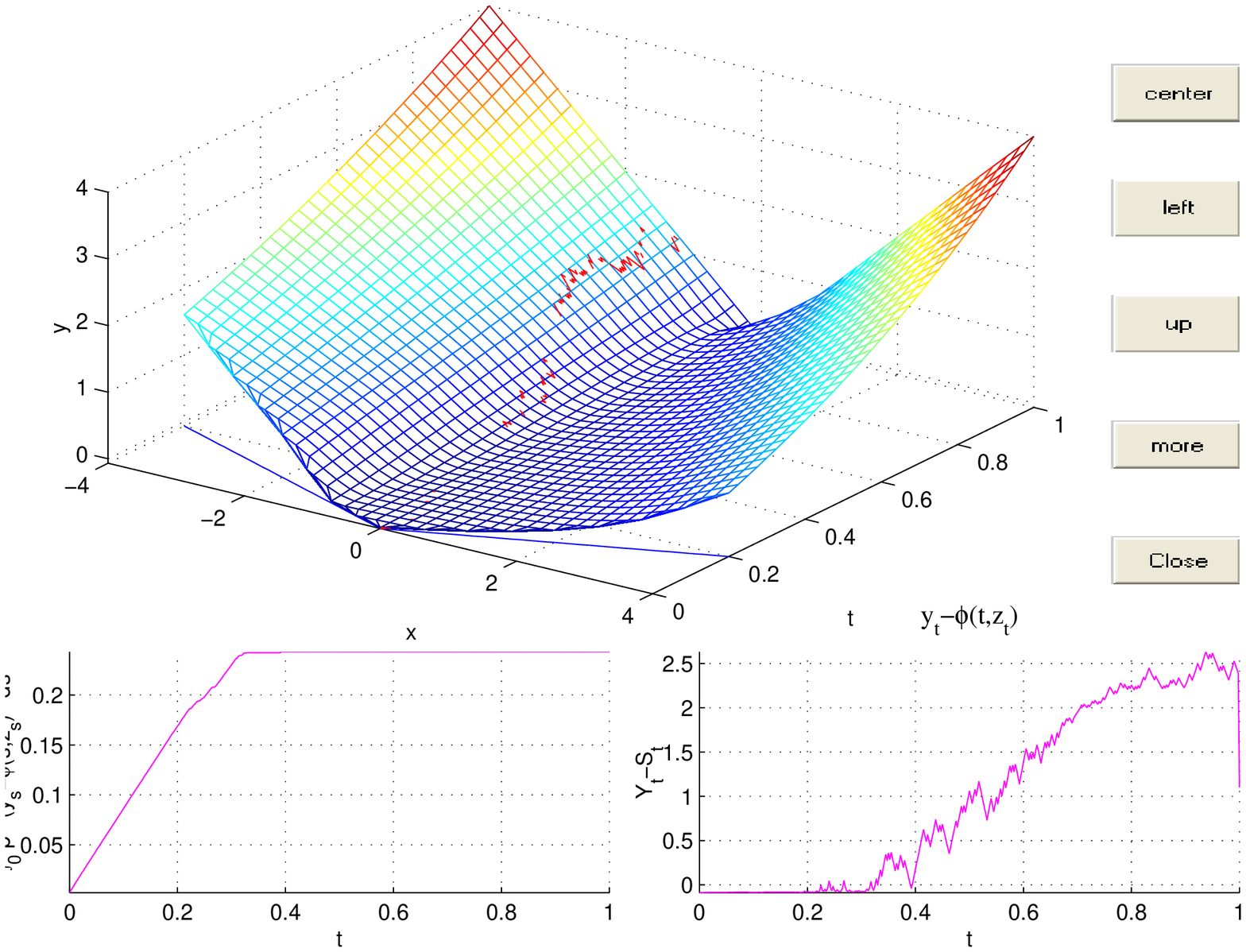} \\[0pt]
Figure 7: The solution surface of penalization BSDE with $d_{\Gamma
}(y,z)=(y-\phi (z))^{-}.$

\centering\includegraphics[height=100mm,width=160mm]{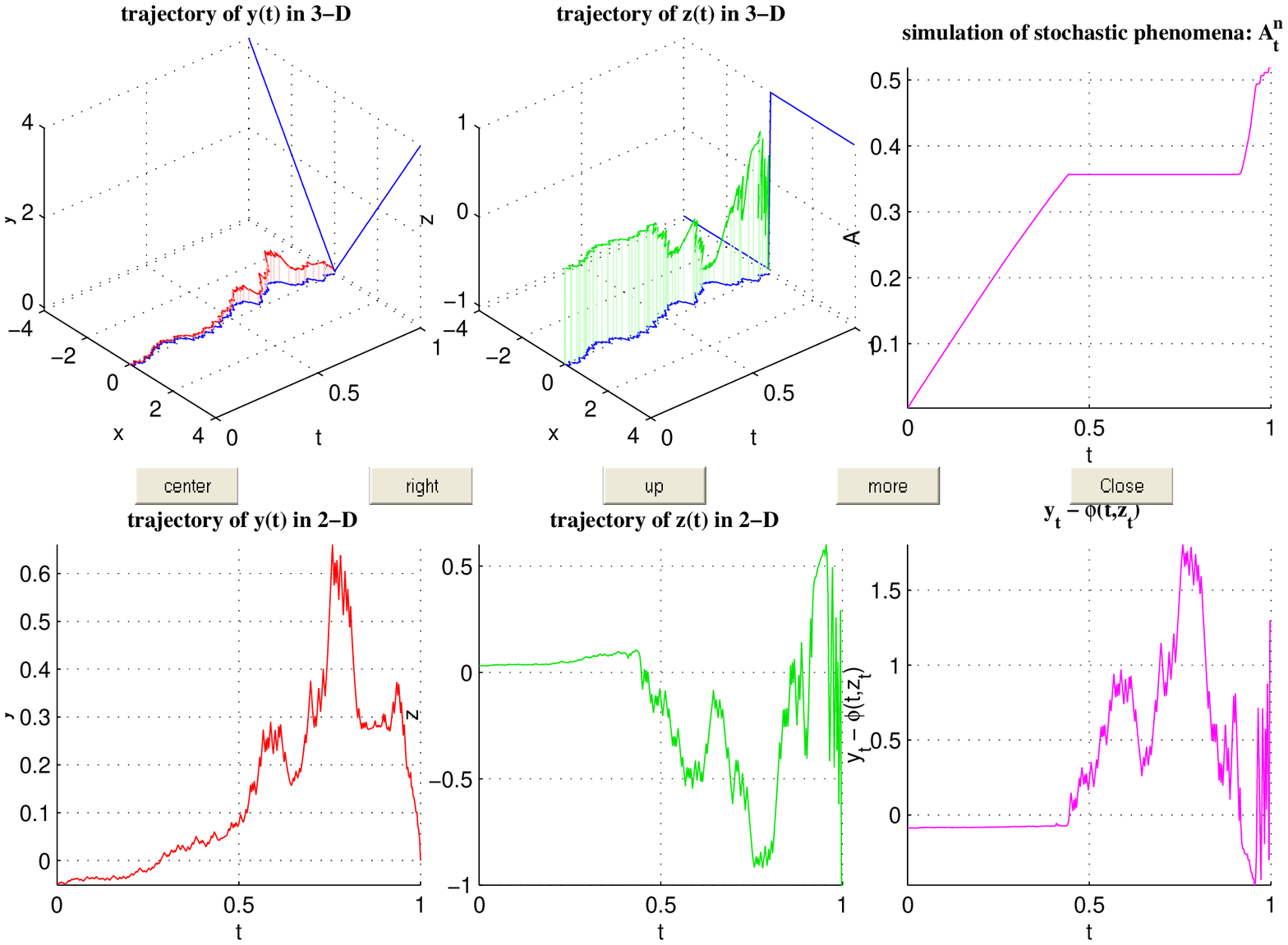} \\[0pt]
Figure 8: A trajectory of solutions of penalization BSDE with
$d_{\Gamma }(y,z)=(y-\phi (z))^{-}.$
\end{center}

\textbf{Acknowledgements: }We appreciate the anonymous referees for
valuable remarks and suggestions.

\end{document}